\documentclass[preprint,12pt,authoryear]{elsarticle}

\usepackage{amssymb,amsmath,amsthm}
\usepackage{bm}
\usepackage{mathtools}
\usepackage{algorithmic}
\usepackage{algorithm}
\usepackage{url}
\usepackage{booktabs}
\usepackage[caption=false]{subfig}
\usepackage{nidanfloat}

\theoremstyle{definition}
\newtheorem{thm}{Theorem}
\newtheorem*{thm*}{Theorem}

\newcommand{\1}{\mbox{1}\hspace{-0.25em}\mbox{l}}

\makeatletter
\def\ps@pprintTitle{%
  \let\@oddhead\@empty
  \let\@evenhead\@empty
  \def\@oddfoot{}%
  \let\@evenfoot\@oddfoot}
\makeatother
\begin{document}

\begin{frontmatter}

\title{Robust personalized pricing under uncertainty of purchase probabilities}

\author[1,3]{Shunnosuke Ikeda}
\author[1]{Naoki Nishimura}

\author[2]{Noriyoshi Sukegawa}

\author[4]{Yuichi Takano}

\affiliation[1]{organization={Product Development Management Office, Data Management \& Planning Office, Recruit Co., Ltd.},
    addressline={1--9--2}, 
    city={Chiyoda--ku},
    % citysep={}, % Uncomment if no comma needed between city and postcode
    postcode={100--6640}, 
    state={Tokyo},
    country={Japan}}
\affiliation[2]{organization={Department of Advanced Sciences, Hosei University},
    addressline={3--7--2 Kajinocho}, 
    city={Koganei--shi},
    % citysep={}, % Uncomment if no comma needed between city and postcode
    postcode={184--8584}, 
    state={Tokyo},
    country={Japan}}
\affiliation[3]{organization={Graduate School of Science and Technology, University of Tsukuba},
    addressline={1--1--1 Tennodai}, 
    city={Tsukuba--shi},
    % citysep={}, % Uncomment if no comma needed between city and postcode
    postcode={305--8573}, 
    state={Ibaraki},
    country={Japan}}
\affiliation[4]{organization={Institute of Systems and Information Engineering, University of Tsukuba},
    addressline={1--1--1 Tennodai}, 
    city={Tsukuba--shi},
    % citysep={}, % Uncomment if no comma needed between city and postcode
    postcode={305--8573}, 
    state={Ibaraki},
    country={Japan}}

\begin{abstract}
This paper is concerned with personalized pricing models aimed at maximizing the expected revenues or profits for a single item.
While it is essential for personalized pricing to predict the purchase probabilities for each consumer, these predicted values are inherently subject to unavoidable errors that can negatively impact the realized revenues and profits.
To address this issue, we focus on robust optimization techniques that yield reliable solutions to optimization problems under uncertainty.
Specifically, we propose a robust optimization model for personalized pricing that accounts for the uncertainty of predicted purchase probabilities.
This model can be formulated as a mixed-integer linear optimization problem, which can be solved exactly using mathematical optimization solvers.
We also develop a Lagrangian decomposition algorithm combined with line search to efficiently find high-quality solutions for large-scale optimization problems.
Experimental results demonstrate the effectiveness of our robust optimization model and highlight the utility of our Lagrangian decomposition algorithm in terms of both computational efficiency and solution quality.
\end{abstract}


% \begin{highlights}
% \item 
% Robust optimization model for personalized pricing considering the uncertainty of predicted purchase probabilities.
% \item MILO formulation of the resultant optimization model.
% \item Scalable heuristic algorithm based on the Lagrangian relaxation method and the golden section search.
% \item Effectiveness of our framework evaluated through simulation experiments.
% \end{highlights}


\end{frontmatter}


\section{Introduction}
\subsection{Background}
In recent years, the rapid advancement of information and communication technologies has fostered a business environment conducive to the collection and retention of vast amounts of data.
Within this context, the significance of utilizing consumer and item information has grown substantially, leading to increased interest in personalized pricing strategies \citep{elmachtoub2021value}.
If the purchasing intentions of each consumer can be accurately predicted, it is possible to maximize revenues and profits by setting personalized prices based on these predictions.

% In response to this trend, numerous industries have conducted experiments with personalized pricing, including grocery chains \citep{clifford2012shopper}, department stores \citep{d2017neiman}, and airlines \citep{tuttle2013flight}, among others \citep{house2015big}.
In response to this trend, numerous industries have conducted experiments with personalized pricing, including grocery chains, department stores, and airlines, among others \citep{elmachtoub2021value,golrezaei2014real}.
Furthermore, personalized pricing has already become an established practice in various sectors, such as professional graduate programs, financial and insurance items, services, and information goods \citep{feldman2015pricing,waldfogel2015first}.

The framework of personalized pricing consists of two phases: demand forecasting and price optimization.
The demand forecasting phase involves constructing a model to predict the probabilities that each consumer will purchase an item at various price points.
Subsequently, the price optimization phase entails solving an optimization problem to determine the price that maximizes a given objective function (e.g., total revenue or profit) based on the constructed prediction model for each consumer.

However, the predicted purchase probabilities inherently include prediction errors, and the uncertainty associated with these errors can significantly reduce the effect of price optimization \citep{duchi2021statistics,delage2010distributionally}.
Consequently, improving the robustness of optimal prices against prediction errors remains a critical challenge in the implementation of personalized pricing strategies.

\subsection{Related work}
Numerous prior studies have been conducted on personalized pricing for a single item \citep{rossi1996value,allenby1998marketing,den2014simultaneously,kallus2021fairness,elreedy2021novel}.
While many of these studies have assumed parametric models for the demand function, it is desirable to use more expressive models without such shape restrictions because the accuracy of demand forecasting directly impacts the effect of price optimization.

Several studies have employed demand forecasting models that do not assume a specific functional form \citep{elmachtoub2021value,biggs2021model,subramanian2022constrained,amram2022optimal,chen2023personalized}.
\citet{elmachtoub2021value} provided lower and upper bounds on the ratio between the profits of an idealized personalized pricing strategy and a single price strategy.
\citet{biggs2021model}, \citet{subramanian2022constrained}, and \citet{amram2022optimal} proposed prescriptive trees to improve the interpretability of optimal prices.
Furthermore, \citet{chen2023personalized} introduced a stochastic optimization model focusing on the fairness of price setting.

\citet{thiele2006single} studied demand uncertainty in the pricing of a single item.
However, this research focused on demand uncertainty associated with time series rather than in the context of personalized pricing.
Robust optimization methods have been proposed for handling demand uncertainty for personalized pricing of multiple items \citep{chen2023model,jin2021price}.
These studies focused on assortment-specific problem structures and had difficulty coping with constraints for a large number of consumers.
% These studies, however, focused on assortment-specific problem structures and are limited to scenarios without consumer-independent constraints.
% Specifically, they do not account for constraints that apply universally across all consumers.

A more general framework for prescriptive analytics based on distributionally robust optimization methods has been proposed by \citet{bertsimas2022bootstrap}, and \citet{wang2021distributionally}.
However, these methods are not readily applicable to price optimization problems where purchase probabilities depend on both consumer covariates and item prices (i.e., decision variable).
Consequently, there is a need for a more comprehensive framework to cope with demand uncertainty in the context of personalized pricing for a single item.

\subsection{Our contribution}
The goal of this paper is to enhance the performance of personalized pricing by employing robust optimization techniques.
Specifically, we propose a robust optimization model for personalized pricing that accounts for the uncertainty of predicted purchase probabilities, where the uncertainty set of purchase probabilities is estimated using calculated using the bootstrap method.
Our optimization model is formulated as a mixed-integer linear optimization problem, which can be solved exactly using mathematical optimization solvers.
Additionally, we develop a Lagrangian decomposition algorithm combined with line search to efficiently find high-quality solutions for large-scale optimization problems.

To validate the effectiveness of our method for robust personalized pricing, we conduct numerical experiments using both synthetic and real-world datasets.
The experimental results demonstrate that our method can significantly enhance expected revenues compared to conventional approaches that do not consider uncertainty.
These findings suggest that more reliable and profitable pricing strategies can be developed using our robust optimization model.
Furthermore, we highlight the utility of our Lagrangian decomposition algorithm in terms of its computational efficiency and solution quality.
The experimental results show that our algorithm can provide near-optimal solutions in considerably less time than that required for exact solutions, enabling the application of our robust optimization method to large-scale price optimization problems.

\section{Personalized pricing model}
In this section, we provide a concise overview of the existing framework for personalized pricing of a single item.
The framework consists of two phases: (i) constructing a model to predict each consumer's purchase probabilities for an item, and (ii) solving a price optimization problem based on purchase probabilities predicted by the constructed model.

Let $p_i$ denote the price of an item offered to consumer $i \in \mathcal{I}$, and ${\bm x}_i$ represent the covariate vector of consumer $i$.
The function $q({\bm x}_i, p_i)$ represents the purchase probability, taking the covariates ${\bm x}_i$ and the offered price $p_i$ as arguments.
We assume that the price offered to one consumer does not affect the purchase probability of other consumers, which reduces the complexity of the optimization problem.
Due to the presence of unique covariate vectors ${\bm x}_i$ for a large number of consumers $i \in \mathcal{I}$, nonlinear and highly accurate prediction models, such as gradient boosting decision trees, are often employed to predict the purchase probabilities.

We define the price vector as ${\bm p} \coloneqq (p_i)_{i \in \mathcal{I}} \in \mathbb{R}^{\lvert \mathcal{I} \rvert}$ and denote the feasible region of the price vector ${\bm p}$ by $\mathcal{P} \subseteq \mathbb{R}^{\lvert \mathcal{I}\rvert}$.
The optimization model for personalized pricing can generally be formulated as follows:
\begin{align}
\underset{{\bm p}}{\parbox{5em}{maximize}} & \sum_{i \in \mathcal{I}} p_i q({\bm x}_i, p_i) \label{eq:obj_per} \\
\parbox{5em}{subject to} & {\bm p} \in \mathcal{P}, \label{eq:price constraints per}
\end{align}
where the objective function in Eq.~\eqref{eq:obj_per} represents the expected total revenue earned from all consumers.

Assuming price $p_i$ is chosen from a set of price candidates $\{ P_{ij} \mid j \in \mathcal{J} \}$, we define a 0--1 decision variable vector ${\bm z} \coloneqq (z_{ij})_{(i,j) \in \mathcal{I} \times \mathcal{J}} \in \{0,1\}^{\lvert \mathcal{I} \times \mathcal{J} \rvert}$ to represent the choice of prices presented to consumers as follows:
\[
z_{ij} = 1 \iff p_{i} = P_{ij}, \quad \forall i \in \mathcal{I}, ~\forall j \in \mathcal{J}.
\]
We also denote the types of constraints by $k \in \mathcal{K}$, and introduce a constraint function $f_{k}({\bm z})$ for $k \in \mathcal{K}$ to represent the feasible region $\mathcal{P}$ of prices.
Consequently, the optimization problem defined in Eqs.~\eqref{eq:obj_per} and \eqref{eq:price constraints per} can be reformulated as follows:
\begin{align}
\underset{{\bm z}~~~}{\parbox{5em}{maximize}} & \sum_{i \in \mathcal{I}}\sum_{k \in \mathcal{K}} P_{ij} q({\bm x}_i, P_{ij})z_{ij} \label{eq:obj_per_dis} \\
\parbox{5em}{subject to} & \sum_{j \in \mathcal{J}} z_{ij} = 1, \quad \forall i \in \mathcal{I}, \label{eq:z_ij=1}\\
& f_{k}({\bm z}) \leq 0, \quad \forall k \in \mathcal{K}, \label{eq:f_k}\\
& z_{ij} \in \{0, 1\}, \quad {\forall} i \in \mathcal{I}, ~\forall j \in \mathcal{J}. \label{eq:z_ij}
\end{align}
The objective function in Eq.~\eqref{eq:obj_per_dis} represents the expected total revenue, similar to Eq.~\eqref{eq:obj_per}.
Eq.~\eqref{eq:z_ij=1} indicates that one price is selected from its candidate set.

\section{Robust optimization model for personalized pricing}
In this section, we present a robust optimization model designed for personalized pricing that effectively addresses the uncertainty associated with the purchase probabilities.
Additionally, we describe the procedure of estimating the uncertainty of the purchase probabilities using the bootstrap method.

\subsection{Problem formulation}
For the consumer $i \in \mathcal{I}$ offered the price $P_{ij}$, let $q_{ij}$ denote the purchase probability with uncertainty, and $\hat{q}_{ij}$ represent the predicted value of the purchase probability.
Let $\Delta_{ij}$ also denote the magnitude of the uncertainty in the purchase probability.
Additionally, $\gamma_i$ is a decision variable indicating the extent to which uncertainty is considered for consumer $i$, and $\Gamma$ represents the maximum number of consumers affected by uncertainty.

Employing the robust optimization technique \citep{bertsimas2004price}, we can formulate the robust optimization model for personalized pricing that accounts for the uncertainty of the purchase probabilities as follows:
\begin{align}
 \underset{{\bm z}~~~}{\parbox{5em}{maximize}}\hspace{-2mm}\underset{{\bm q},{\bm \gamma}~~~}{\parbox{5em}{minimize}} & \sum_{i \in \mathcal{I}}\sum_{j \in \mathcal{J}} P_{ij} q_{ij} z_{ij} \label{eq:obj_robust}\\
 \parbox{5em}{subject to} & \text{Eqs.}~\eqref{eq:z_ij=1}\text{--}\eqref{eq:z_ij}, \\
 & q_{ij} = \hat{q}_{ij} - \gamma_i \Delta_{ij}, \quad \forall i \in \mathcal{I}, ~\forall j \in \mathcal{J}, \label{eq:q_hat} \\ 
 & 0 \leq \gamma_i \leq 1, \quad \forall i \in \mathcal{I}, \label{eq:gamma}\\
 & \sum_{i \in \mathcal{I}}\gamma_i \leq \Gamma, \label{eq:sum_gamma}
\end{align}
where ${\bm q} \coloneqq (q_{ij})_{(i, j) \in \mathcal{I} \times \mathcal{J}}$ and ${\bm \gamma} \coloneqq (\gamma_i)_{i \in \mathcal{I}}$.
Eq.~\eqref{eq:obj_robust} seeks to maximize the expected revenue in the worst case.
Eq.~\eqref{eq:q_hat} represents the purchase probability taking uncertainty into account.
Since the purchase probability $q_{ij}$ is non-negative, it is necessary that $0 \leq \Delta_{ij} \leq \hat{q}_{ij}$ to ensure feasibility.
Eq.~\eqref{eq:gamma} indicates the degree of impact of the uncertainty for each consumer, while Eq.~\eqref{eq:sum_gamma} limits the number of consumers affected by the uncertainty.

However, the optimization problem~\eqref{eq:obj_robust}--\eqref{eq:sum_gamma} is difficult to solve directly using mathematical optimization solvers because it requires maximizing the minimum value of the objective function.
Therefore, we reformulate it into a more tractable form by converting the inner minimization problem into its dual problem, allowing for an easier solution.

Now, let $\mu_{i} \in \mathbb{R}_{+}$ and $\nu \in \mathbb{R}_{+}$ denote the Lagrange multipliers corresponding to Eqs.~\eqref{eq:gamma} and \eqref{eq:sum_gamma}, respectively.
The following theorem is then established:
\begin{thm} \label{thm:robust}
The optimization problem~\eqref{eq:obj_robust}--\eqref{eq:sum_gamma} can be reformulated as the following MILO problem~\eqref{eq:obj_robust_dual}--\eqref{eq:nu}.
\end{thm}
\vspace{-8mm}
\begin{align}
 \underset{{\bm z},{\bm \mu},\nu~~}{\parbox{5em}{maximize}} & \sum_{i \in \mathcal{I}} \biggl(\sum_{j \in \mathcal{J}}P_{ij} \hat{q}_{ij} z_{ij} - \mu_{i}\biggr) - \Gamma \nu \label{eq:obj_robust_dual}\\
 \parbox{5em}{subject to} & \text{Eqs.}~\eqref{eq:z_ij=1}\text{--}\eqref{eq:z_ij}, \\
 & \mu_i + \nu - \sum_{j \in \mathcal{J}}P_{ij} \Delta_{ij} z_{ij} \geq 0, \quad \forall i \in \mathcal{I}, \label{eq:mu_nu}\\ 
 & \mu_i \geq 0, \quad \forall i \in \mathcal{I}, \label{eq:mu}\\
 & \nu \geq 0,  \label{eq:nu}
\end{align}
where ${\bm \mu} \coloneqq (\mu_i)_{i \in \mathcal{I}}$.
\begin{proof}
Substituting the constraint~\eqref{eq:q_hat} into the objective function~\eqref{eq:obj_robust}, the inner minimization problem can be written as follows:
\begin{align}
\underset{{\bm \gamma}~~~}{\parbox{5em}{minimize}} & \sum_{i \in \mathcal{I}}\sum_{j \in \mathcal{J}} \biggl(P_{ij} \hat{q}_{ij} z_{ij} - P_{ij} \gamma_i \Delta_{ij} z_{ij}\biggr) \label{eq:obj_inner} \\
 \parbox{5em}{subject to} & \text{Eqs.}~\eqref{eq:gamma} ~\text{and}~\eqref{eq:sum_gamma}. \label{eq:10_11}
\end{align}
The dual problem of this inner minimization problem can be written as follows:
\begin{align*}
\underset{{\bm \mu}, \nu~~}{\parbox{5em}{maximize}} & \sum_{i \in \mathcal{I}}\biggl(\sum_{j \in \mathcal{J}} P_{ij} \hat{q}_{ij} z_{ij} - \mu_{i}\biggr) - \Gamma \nu \\
\parbox{5em}{subject to} & \mu_{i} + \nu - \sum_{j \in \mathcal{J}} P_{ij} \Delta_{ij} z_{ij} \geq 0, \quad \forall i \in \mathcal{I}, \\
& \mu_{i} \geq 0, \quad \forall i \in \mathcal{I}, \\
& \nu \geq 0.
\end{align*}
Therefore, the optimization problem \eqref{eq:obj_robust}--\eqref{eq:sum_gamma} can be transformed into the maximization problem \eqref{eq:obj_robust_dual}--\eqref{eq:nu}, completing the proof.
\end{proof}

According to Theorem \ref{thm:robust}, the optimization problem~\eqref{eq:obj_robust}--\eqref{eq:sum_gamma} can be converted to a MILO problem, which can be solved exactly using optimization solvers.

\subsection{Uncertainty estimation of purchase probabilities} 
We here describe the procedure for calculating the parameter $\Delta_{ij}$, which indicates the magnitude of uncertainty in the purchase probability when the price candidate $P_{ij}$ is offered to the consumer $i$, using the bootstrap method \citep{efron1992bootstrap}.

Specifically, we construct multiple purchase probability prediction models using training data replicated through the bootstrap method.
Then, $\Delta_{ij}$ is calculated based on the standard deviation of the predicted values obtained from these prediction models.
This approach enables the assessment of the magnitude of uncertainty associated with any prediction model, providing a robust foundation for our personalized pricing strategy.

Let $\mathcal{I}^{{\rm train}}$ denote the set of consumers in the training data used to construct the prediction model of purchase probabilities.
Let $N^{{\rm bs}}$ represent the number of trials in the bootstrap method, and let $\hat{q}^{(n^{{\rm bs}})}_{ij}$ represent the predicted purchase probability for consumer $i$ if price $P_{ij}$ is offered in each trial $n^{{\rm bs}} \in \{1, 2, \ldots, N^{{\rm bs}}\}$.
The standard deviation of the predicted purchase probabilities is denoted by $\hat{\sigma}_{ij}$.

Furthermore, $\kappa \in \mathbb{R}_{+}$ is a parameter that indicates the ratio to the calculated standard deviation $\hat{\sigma}_{ij}$, and it defines the magnitude of uncertainty $\Delta_{ij}$ in the purchase probability.
For example, when $\kappa = 1$, the magnitude of uncertainty $\Delta_{ij}$ in the purchase probability is equal to the estimated standard deviation.

The detailed procedure for calculating the parameter $\Delta_{ij}$ is presented below.
\vspace{3mm}

\noindent
\textbf{Calculation procedure for $\Delta_{ij}$}
\begin{itemize}
    \setlength{\leftskip}{0.7cm}
    \setlength{\itemsep}{5pt}
    \item[\textrm{Step 1:}] Extract $\lvert \mathcal{I}^{{\rm train}} \rvert$ data instances with replacement from the training data $\mathcal{I}^{{\rm train}}$ to generate a bootstrap sample. This process is repeated $N^{{\rm bs}}$ times to produce $N^{{\rm bs}}$ bootstrap samples.
    \item[\textrm{Step 2.}] Construct prediction models for each of the $N^{{\rm bs}}$ bootstrap samples generated in Step 1.
    \item[\textrm{Step 3.}] Calculate the standard deviation of the predicted purchase probabilities for the target consumer $i \in \mathcal{I}^{\textrm{test}}$ and the price candidate index $j \in \mathcal{J}$ as follows:
    \begin{align*}
        &\bar{q}_{ij} = \frac{1}{N^{{\rm bs}}}\sum_{n^{{\rm bs}}=1}^{N^{{\rm bs}}} \hat{q}^{(n^{{\rm bs}})}_{ij}, \quad \forall i \in \mathcal{I}^{\textrm{test}}, ~ \forall j \in \mathcal{J}, \\
        &\hat{\sigma}_{ij} = \sqrt{\frac{1}{N^{{\rm bs}}-1}\sum_{n^{{\rm bs}}=1}^{N^{{\rm bs}}} \biggl(\hat{q}^{(n^{{\rm bs}})}_{ij} - \bar{q}_{ij})\biggr)^2}, \quad \forall i \in \mathcal{I}^{\textrm{test}}, ~ \forall j \in \mathcal{J}.
    \end{align*}
    \item[\textrm{Step 4.}] Calculate $\Delta_{ij}$ based on the following equation. The minimum value with is selected to ensure that $q_{ij} \geq 0$.
    \begin{align*}
        \Delta_{ij} =  \min(\kappa \hat{\sigma}_{ij}, \hat{q}_{ij}), \quad \forall i \in \mathcal{I}^{\textrm{test}}, ~ \forall j \in \mathcal{J}.
    \end{align*}
\end{itemize}

\section{Heuristic algorithm for large-scale problems}
In this section, we design a scalable heuristic algorithm for efficiently finding high-quality solutions to a large-scale problem, combining the Lagrangian relaxation method \citep{geoffrion2009lagrangean, guignard2003lagrangean} and the golden section search \citep{kiefer1953sequential}.

By decomposing the original problem~\eqref{eq:obj_robust_dual}--\eqref{eq:nu} into small optimization problems for consumers, our algorithm achieves fast and high-quality solutions even as the number of consumers increases.
The original problem is structured to be indivisible due to the constraints~\eqref{eq:f_k} and the decision variable $\nu$, which depend on the decision variables $\bm{z}$.

\subsection{Lagrangian relaxation method} \label{sec:lagrange_relaxation_method}
First, to address the constraints~\eqref{eq:f_k}, we introduce Lagrangian multipliers as ${\bm \lambda} \coloneqq (\lambda_k)_{k \in \mathcal{K}} \in \mathbb{R}^{\lvert \mathcal{K} \rvert}_{+}$ and consider a Lagrangian relaxation problem by relaxing the constraints~\eqref{eq:f_k} as follows:
\begin{align}
 \underset{{\bm z},{\bm \mu},{\bm \nu}~~~}{\parbox{5em}{maximize}} & \sum_{i \in \mathcal{I}} \biggl(\sum_{j \in \mathcal{J}}P_{ij} \hat{q}_{ij} z_{ij} - \mu_{i}\biggr) - \Gamma \nu - \sum_{k \in \mathcal{K}} \lambda_{k} f_k({\bm z})\label{eq:obj_lagrange}\\
 \parbox{5em}{subject to} & \textrm{Eqs.}~\eqref{eq:z_ij=1}, \eqref{eq:z_ij}, \eqref{eq:mu_nu}\text{--}\eqref{eq:nu}. \label{eq:cons_lagrange} 
\end{align}
The optimal value of the Lagrangian relaxation problem serves as an upper bound for the optimal value of the original problem.

To minimize the upper bound, we consider the following Lagrange dual problem.
\begin{align}
 \underset{{\bm \lambda}~~~}{\parbox{5em}{minimize}}\hspace{-2mm}\underset{{\bm z},{\bm \mu},{\bm \nu}~~~}{\parbox{5em}{maximize}} & \sum_{i \in \mathcal{I}} \biggl(\sum_{j \in \mathcal{J}}P_{ij} \hat{q}_{ij} z_{ij} - \mu_{i}\biggr) - \Gamma \nu - \sum_{k \in \mathcal{K}} \lambda_{k} f_k({\bm z})\label{eq:obj_lagrange_dual}\\
 \text{subject to}~~ & \textrm{Eqs.}~\eqref{eq:z_ij=1}, \eqref{eq:z_ij}, \eqref{eq:mu_nu}\text{--}\eqref{eq:nu} \label{eq:cons_lagrange_dual} 
\end{align}
Due to the structure of the Lagrangian dual problem, which embeds a maximization problem within a minimization problem, it is difficult to directly solve the problem.

Therefore, we employ the projected subgradient method~\citep{boyd2003subgradient}, to update the Lagrange multipliers ${\bm \lambda}$.
With the updated multipliers, we then solve the Lagrangian relaxation problem.
By alternately repeating this process until the termination conditions are satisfied, it becomes possible to approximately solve the Lagrangian dual problem.
% \vspace{3mm}
% \noindent
% \textbf{Updating ${\bm \lambda}$ using the projected subgradient method}
% \vspace{3mm}

At the $t$-th iteration, let the step size parameter be $\delta^{(t)}$, the solutions be ${\bm z}^{(t)}$, and the Lagrange multipliers be ${\bm \lambda}^{(t)}$.
By denoting the objective function of Eq.~\eqref{eq:obj_lagrange} as $\rho({\bm \lambda}^{(t)})$, the update formula for the Lagrange multipliers can be described as follows:
\begin{align}
    \lambda^{(t+1)}_{k} &\leftarrow \max \biggl(\lambda^{(t)}_{k} - \delta^{(t)}\frac{\partial}{\partial \lambda^{(t)}_k}\rho({\bm \lambda}^{(t)}), 0 \biggr), \quad \forall k \in \mathcal{K} \notag \\
    \iff \lambda^{(t+1)}_{k}& \leftarrow \max \biggl(\lambda^{(t)}_{k} + \delta^{(t)}f_k({\bm z}^{(t)}), 0 \biggr), \quad \forall k \in \mathcal{K}. \label{eq:lagrange_update}
\end{align}
Thus, by using the Lagrangian relaxation method, it becomes possible to incorporate the constraint~\eqref{eq:f_k} into the objective function, thereby reducing the computational complexity associated with increasing the number of the decision variables $\bm{z}$.

\subsection{Golden section search} \label{sec:golden_section_search}
Second, to address the decision variable $\nu$, we employ the golden section search.
By fixing the decision variable $\nu$ as a constant, the Lagrangian relaxation problem (Eqs.~\eqref{eq:obj_lagrange} and \eqref{eq:cons_lagrange}) can be decomposed into optimization problems for each consumer $i$ as follows:
\begin{align}
 \underset{{\bm z}_{i},\mu_{i}~~~}{\parbox{5em}{maximize}} & \sum_{j \in \mathcal{J}}P_{ij} \hat{q}_{ij} z_{ij} - \mu_{i} - \sum_{k \in \mathcal{K}} \lambda_{k} f_k({\bm z}_{i}) - \frac{\Gamma \nu}{\lvert\mathcal{I}\rvert} \label{eq:obj_lagrange_split} \\
 \parbox{5em}{subject to} & \mu_i + \nu - \sum_{j \in \mathcal{J}}P_{ij} \Delta_{ij} z_{ij} \geq 0,  \\
 & \mu_i \geq 0, \\
 & \sum_{j \in \mathcal{J}}z_{ij} = 1, \\
 & z_{ij} \in \{0, 1\}, \quad \forall j \in \mathcal{J}, \label{eq:cons_lagrange_split}
\end{align}
where ${\bm z}_i \coloneqq (z_{ij})_{j \in \mathcal{J}} \in \{0, 1\}^{\vert\mathcal{J}\rvert}$.
This optimization problem can be solved to select prices from the set of candidate prices for each consumer $i$ so as to maximize the objective function while satisfying the constraints.
This approach allows to obtain solutions with a computational complexity of $\mathcal{O}(\lvert \mathcal{I} \rvert \lvert \mathcal{J} \rvert)$.
Since these optimization problems can be solved independently for each consumer $i$, it is possible to accelerate the process through parallel processing.

Given that the constant $\nu$ is a one-dimensional value, we adopt the golden section search (Algorithm \ref{alg1}) to iteratively solve the optimization problem~\eqref{eq:obj_lagrange_split}--\eqref{eq:cons_lagrange_split} and thus identify the appropriate value for $\nu$.
For a given $\nu$, let the optimal value of Eq.~\eqref{eq:obj_lagrange_split} denote $\pi^{*}_{i}(\nu)$, and we define $\pi^{*}(\nu) \coloneqq \sum_{i \in \mathcal{I}} \pi^{*}_{i}(\nu)$.
\begin{algorithm}[H]
\caption{Golden Section Search}
\label{alg1}
\begin{algorithmic}[1]
    \STATE Set the lower and upper bounds of the search interval $[a, b]$.
    \STATE Set the golden ratio $\phi = \frac{1 + \sqrt{5}}{2}$.
    \STATE Set $c \leftarrow b - \frac{b - a}{\phi}$.
    \STATE Set $d \leftarrow a + \frac{b - a}{\phi}$.
    \STATE Repeat the following steps until $\lvert a - b \rvert < \varepsilon^{g}$ for a given small positive constant $\varepsilon^{g}$:
        \STATE \hspace{\algorithmicindent} If $\pi^{*}(c) < \pi^{*}(d)$, then set $b \leftarrow d$.
        \STATE \hspace{\algorithmicindent} Otherwise, set $a \leftarrow c$.
        \STATE \hspace{\algorithmicindent} Update $c \leftarrow b - \frac{b - a}{\phi}$.
        \STATE \hspace{\algorithmicindent} Update $d \leftarrow a + \frac{b - a}{\phi}$.
    \STATE Output $\frac{a+b}{2}$ as the solution.
\end{algorithmic}
\end{algorithm}

\subsection{Algorithm procedure}
By integrating the Lagrangian relaxation method (Section \ref{sec:lagrange_relaxation_method}) with the golden section search (Section \ref{sec:golden_section_search}), the original problem~\eqref{eq:obj_robust_dual}--\eqref{eq:nu} can be decomposed into optimization problems for each consumer $i$, allowing to obtain fast and of high-quality solutions.

Let $\nu^{(t)}$ denote $\nu$ at iteration $t$, and Algorithm \ref{alg2} gives a detailed procedure of our proposed algorithm.
Our algorithm first initializes the Lagrangian multipliers ${\bm \lambda}^{(t)}$ and the step size $\delta^{(t)}$ (lines 2--3).
It then calculates $\nu^{(t)}$ using the golden section search (line 5).
Subsequently, the algorithm computes the subgradient by solving the optimization problem~\eqref{eq:obj_lagrange_split}--\eqref{eq:cons_lagrange_split} given ${\bm \lambda}^{(t)}, \nu^{(t)}$ and update Lagrange multipliers ${\bm \lambda}^{(t)}$ and the step size $\delta^{(t)}$ (lines 6--8).
These steps are repeated until the termination conditions are satisfied (lines 4--9).

Recall that our algorithm can be accelerated through parallel implementation, as the computations can be performed independently for each consumer $i$.

\begin{algorithm}[H]
\caption{Proposed Algorithm}
\label{alg2}
\begin{algorithmic}[1]
    \STATE Set the iteration count $t = 1$.
    \STATE Initialize the Lagrange multipliers ${\bm \lambda}^{(t)}$.
    \STATE Set the initial step size $\delta^{(t)}$.
    \STATE Repeat the following steps until the termination conditions are satisfied:
        \STATE \hspace{\algorithmicindent} Calculate $\nu^{(t)}$ using the golden section search (Algorithm\ref{alg1}).
        \STATE \hspace{\algorithmicindent} Solve the optimization problems~\eqref{eq:obj_lagrange_split}--\eqref{eq:cons_lagrange_split} given ${\bm \lambda}^{(t)}, \nu^{(t)}$ and compute the subgradient $\nabla_{{\bm \lambda}^{(t)}}\rho({\bm \lambda}^{(t)})$ based on the solution ${\bm z}^{(t)}$.
        \STATE \hspace{\algorithmicindent} Update the Lagrange multipliers ${\bm \lambda}^{(t)}$ according to the formula~\eqref{eq:lagrange_update}.
        \STATE \hspace{\algorithmicindent} Update the step size $\delta^{(t)}$.
        \STATE \hspace{\algorithmicindent} Update $t \leftarrow t + 1$.
    \STATE Output the solution ${\bm z}^{(t)}$.
\end{algorithmic}
\end{algorithm}

\section{Experimental results}
In this section, we evaluate the effectiveness of our proposed method through simulation experiments using synthetic and real-world datasets.
All experiments were performed on a Mac OS 12.6 computer equipped with an Apple M1 chip (8 cores) and 8 GB RAM.
To solve our optimization problems exactly, we used Gurobi Optimizer 10.0.2 \footnote{https://www.gurobi.com/}, a state-of-the-art commercial solver for mathematical optimization.

\subsection{Synthetic datasets}
In the numerical experiments with synthetic datasets, we generated six synthetic datasets from various purchase probability models, following the previous work \citep{biggs2021model}.
By using the true purchase probability model in these experiments, we can accurately calculate counterfactual outcomes following price changes and assess expected revenues.
\vspace{3mm}

\noindent
\textbf{Purchase Probability Model}
\vspace{3mm}

Let the generated price and covariates denote $P$ and ${\bm X} \coloneqq(X_1, X_2, \ldots, X_n)$, respectively.
Furthermore, let $g({\bm X})$ and $h({\bm X})$ represent any transformation functions applied to the covariates ${\bm X}$ and denote the error as $\epsilon$.
The purchase probability $q({\bm X}, P)$ for generating synthetic datasets was thus defined as follows:
\begin{align}
    q({\bm X}, P)^* = g({\bm X}) + h({\bm X}) P + \epsilon, \quad q({\bm X}, P)= \begin{cases}1 & \text { if } q({\bm X}, P)^*>0, \\ 0 & \text { if } q({\bm X}, P)^* \leq 0 .\end{cases}
\end{align}

Table \ref{tab:synthetic_setting} details the purchase probability model for each synthetic dataset.
Let $I_{n} \in \mathbb{R}^{n \times n}$ denote a matrix with diagonal elements of 1 and off-diagonal elements of 0.
The indicator function $\1\{\cdot\}$ returns 1 if the condition within the braces is satisfied, and 0 otherwise.
Moreover, for all datasets, the error term $\epsilon$ was assumed to follow a normal distribution with $\epsilon \sim \mathcal{N}(0, 2)$.

\begin{table}[H]
\centering
\caption{Purchase probability models in synthetic datasets}
\scalebox{0.88}{
  \begin{tabular}{ll}
    \hline
    Dataset & Purchase probability model   \\
    \hline
    Dataset1  & Linear probit model with no confounding ($n=1$):   \\
               & $g({\bm X})=X_1$, $h({\bm X})=-1$, ${\bm X} \sim \mathcal{N}(5, 1)$, $P \sim \mathcal{N}(5,2)$. \\
               \addlinespace[2mm]
    Dataset2  & Higher dimension probit model with sparse linear interaction ($n=20$):  \\
               & $g({\bm X})=5$, $h({\bm X})=-1.5{\bm \beta}{\bm X}^{\top}$, ${\bm X} \sim \mathcal{N}({\bm 0}, I_{20})$, ${\bm \beta} \coloneqq (\beta_1, \ldots, \beta_{20})$\\
               & $(\beta_1, \ldots, \beta_5) \sim \mathcal{N}({\bm 0}, I_{5})$, $\beta_i = 0 ~(i \in \{6, \ldots, 20\})$, $P \sim \mathcal{N}(5,2)$. \\
               \addlinespace[2mm]
    Dataset3  & Probit model with step interaction ($n=1$):  \\
               & $g({\bm X})=5$, $h({\bm X})=-1.2 \1 \{X_1 < -1\} - 1.1\1 \{-1 \leq X_1 < 0\}$  \\
               & $- 0.9\1\{0 \leq X_1 < 1\} - 0.8\1\{1 \leq X_1\}$, ${\bm X} \sim \mathcal{N}(0, 1)$, $P \sim \mathcal{N}(X_1 + 5,2)$.\\
               \addlinespace[2mm]
    Dataset4  & Probit model with multi-dimensional step interaction ($n=2$):\\
               & $g({\bm X})=5$, $h({\bm X})=-1.25 \1 \{X_1 < -1\} - 1.1\1 \{-1 \leq X_1 < 0\} $ \\
               & $- 0.9\1\{0 \leq X_1 < 1\} - 0.75\1\{1 \leq X_1\} + 0.1\1\{X_2 < 0\} - 0.1\1\{X_2 \geq 0\}$ \\
               & ${\bm X} \sim \mathcal{N}({\bm 0}, I_2)$, $P \sim \mathcal{N}(X_1 + 5,2)$. \\
               \addlinespace[2mm]
    Dataset5  & Linear probit model with observed confounding ($n=1$):\\
               & $g({\bm X})=X_1$, $h({\bm X})=-1$, ${\bm X} \sim \mathcal{N}(5, 1)$, $P \sim \mathcal{N}(X_1 + 5,2)$. \\
               \addlinespace[2mm]
    Dataset6  & Probit model with non-linear interaction ($n=2$): \\
               & $g({\bm X}) = 4\lvert X_1 + X_2 \rvert$, $h({\bm X}) = -\lvert X_1 + X_2 \rvert$, \\
               & ${\bm X} \sim \mathcal{N}({\bm 0}, I_2)$, $P \sim \mathcal{N}(X_1 + 5,2)$. \\
    \hline
  \end{tabular}
}
\label{tab:synthetic_setting}
\end{table}

\noindent
\textbf{Evaluation procedure}
\vspace{3mm}

We evaluated the revenues in synthetic datasets according to the following process.
The set of candidate prices was defined as nine different values ranging from the 10th to the 90th percentile of the prices observed in the training data.
\begin{enumerate}
\item Generate $\lvert \mathcal{I}^{{\rm train}} \rvert$ training data samples based on the purchase probability model and train the purchase prediction model;
\item Generate $\lvert \mathcal{I}^{{\rm test}} \rvert$ testing data samples in accordance with the purchase probability model;
\item Input covariates into the purchase prediction model to compute predicted purchase probabilities, and determine prices based on these predictions; 
\item Evaluate revenues for the testing data by calculating the occurrence of purchases according to the purchase probability model associated with the generated covariates and determined prices.
\end{enumerate}

\subsection{Real-world datasets}
In numerical experiments with real-world dataset, as in the previous studies \citep{biggs2021model}, we use purchase data from grocery stores to evaluate our proposed method.
\vspace{3mm}

\noindent
\textbf{Grocery pricing}
\vspace{3mm}

We used publicly available data from Dunnhumby titled ``The complete journey"\footnote{\url{https://www.dunnhumby.com/source-files/}}, which includes two years of transactional grocery purchase data from 2,500 households frequently shopping at grocery stores, along with their demographic information.
In this dataset, available household information includes consumers' discretized age, discretized household income, whether they are a homeowner or renter, and household composition.
Household composition categories are single male, single female, 2 adults no kids, 2 adults with kids, 1 adult with kids, and unknown.

Consistent with the prior work \citep{biggs2021model}, we followed established procedures for preprocessing data and building prediction models of purchase probabilities.
Each record represents a purchase opportunity, detailing household information and the unit price of strawberries ranging from \$1.99 to \$5.00, mostly in \$0.50 increments, along with the outcome of whether strawberries were purchased.
If strawberries were not purchased, the price information is missing and we imputed it using the average of the prices of the last three purchases.

The dataset was divided into 50\% training data and 50\% testing data.
The training data was used to build a prediction model of purchase probabilities, and purchase probabilities were predicted for each price in the testing data to determine optimal pricing.
Since the true purchase probability model was unknown in the real-world dataset, an independent purchase probability model was also constructed for the testing data to evaluate pricing.
LightGBM \citep{ke2017lightgbm} was adopted as the prediction model with the hyperparameter $\mathsf{num\_boost\_round}$ (i.e., the number of boosting iterations) set to 50 and other hyperparameters at default settings.
Candidate prices range from \$1.99 to \$4.99 in \$0.50 increments, for a total of seven prices.

\subsection{Evaluation of robust optimization}
First, we present the experimental results and discussion on both synthetic and real-world datasets to verify the effectiveness of the proposed robust optimization method~\eqref{eq:obj_robust_dual}--\eqref{eq:nu}.
To investigate the effect of varying the number of consumers affected by the uncertainty, we introduced a parameter $\alpha \in [0, 1]$ representing the ratio of consumers affected to all consumers, and set $\Gamma = \alpha\lvert\mathcal{I^{\rm test}}\rvert$.
When $\alpha=0$, the scenario assumes no effect of uncertainty, while $\alpha=1$ corresponds to a situation where all consumers are affected by uncertainty.
Additionally, the number of bootstrap iterations in our method was set to $N^{{\rm bs}} = 20$.

\subsubsection{Comparative method}
\cite{biggs2021model} employed LightGBM as a teacher model and introduced the Greedy Student Prescriptive Tree (SPT), which demonstrated superior performance compared to other competing prescriptive tree methods, such as those proposed by \cite{kallus2017recursive} and \cite{athey2016recursive}.
Building upon this work, \cite{amram2022optimal} proposed Optimal Policy Trees (OPT), enhancing performance by optimizing the branches of the SPT.
\cite{subramanian2022constrained} also extended this approach to a multiway tree that incorporated constraints, achieving performance surpassing that of OPT.

While the performance of these prescriptive tree methods has improved over time, previous approaches have not accounted for uncertainty in predicted values.
Moreover, the most recent method \citep{subramanian2022constrained} has demonstrated performance nearly equivalent to a comparative method called LGBM, which determines prices for consumers based on purchase probabilities predicted by LightGBM.

Our proposed method can be positioned as an extension of LGBM, which has served as a comparative method in previous studies.
Notably, when $\alpha=0$, our method is identical to LGBM.
Therefore, in this numerical experiment, we adopted LGBM (equivalent to our proposed method when $\alpha=0$) as the comparative method, consistent with previous studies.

\subsubsection{Result for synthetic datasets}
We analyze the results of purchase probability predictions and the expected revenue resulting from price optimization.
These results were based on $\lvert \mathcal{I}^{\rm test} \rvert = 500$, and the computation time for price optimization was limited to 600 seconds.
For the purchase probability prediction model, default hyperparameter values were used.
The training process used 20\% of the training data as a validation data and was terminated when no improvement in the area under the curve (AUC) was observed on the validation set for 10 consecutive iterations.

First, we discuss the results of the purchase probability predictions.
Table \ref{tab:auc} shows the AUC of the purchase probability prediction model for each synthetic dataset.
From this table, we observe that as the number of consumers in the training data, $\lvert \mathcal{I}^{{\rm train}} \rvert$, increased, the AUC for the testing data also improved.
However, in the case of Dataset2, the AUC value remained below 0.6, indicating a difficult prediction environment.
\begin{table}[H]
    \centering
    \caption{AUC for testing data in synthetic datasets ($\lvert \mathcal{I}^{\rm{test}} \rvert=500$)}
    \catcode`?=\active \def?{\phantom{0}}% ?:=\phantom{0}
        \begin{tabular}{crr}
        \toprule
          & \multicolumn{2}{c}{$\lvert \mathcal{I}^{{\rm train}} \rvert$} \\
        \cmidrule(lr){2-3} 
        Dataset & $100$ & $1000$  \\
        \midrule
        Dataset1 & $0.784~(\pm0.014)$ & $0.826~(\pm0.004)$ \\
        \addlinespace[2mm]
        Dataset2 & $0.517~(\pm0.015)$ & $0.556~(\pm0.006)$ \\
        \addlinespace[2mm]
        Dataset3 & $0.775~(\pm0.007)$ & $0.810~(\pm0.004)$ \\
        \addlinespace[2mm]
        Dataset4 & $0.781~(\pm0.007)$ & $0.810~(\pm0.006)$ \\
        \addlinespace[2mm]
        Dataset5 & $0.751~(\pm0.008)$ & $0.806~(\pm0.004)$ \\
        \addlinespace[2mm]
        Dataset6 & $0.744~(\pm0.014)$ & $0.804~(\pm0.005)$ \\ \bottomrule
        \end{tabular} 
        \label{tab:auc}
\end{table}

\begin{figure}[H]
\centering
    \subfloat[Dataset1 ($\lvert \mathcal{I}^{{\rm train}}\rvert = 100$)]{\includegraphics[width=2.7in, bb=0 0 607 518,clip]{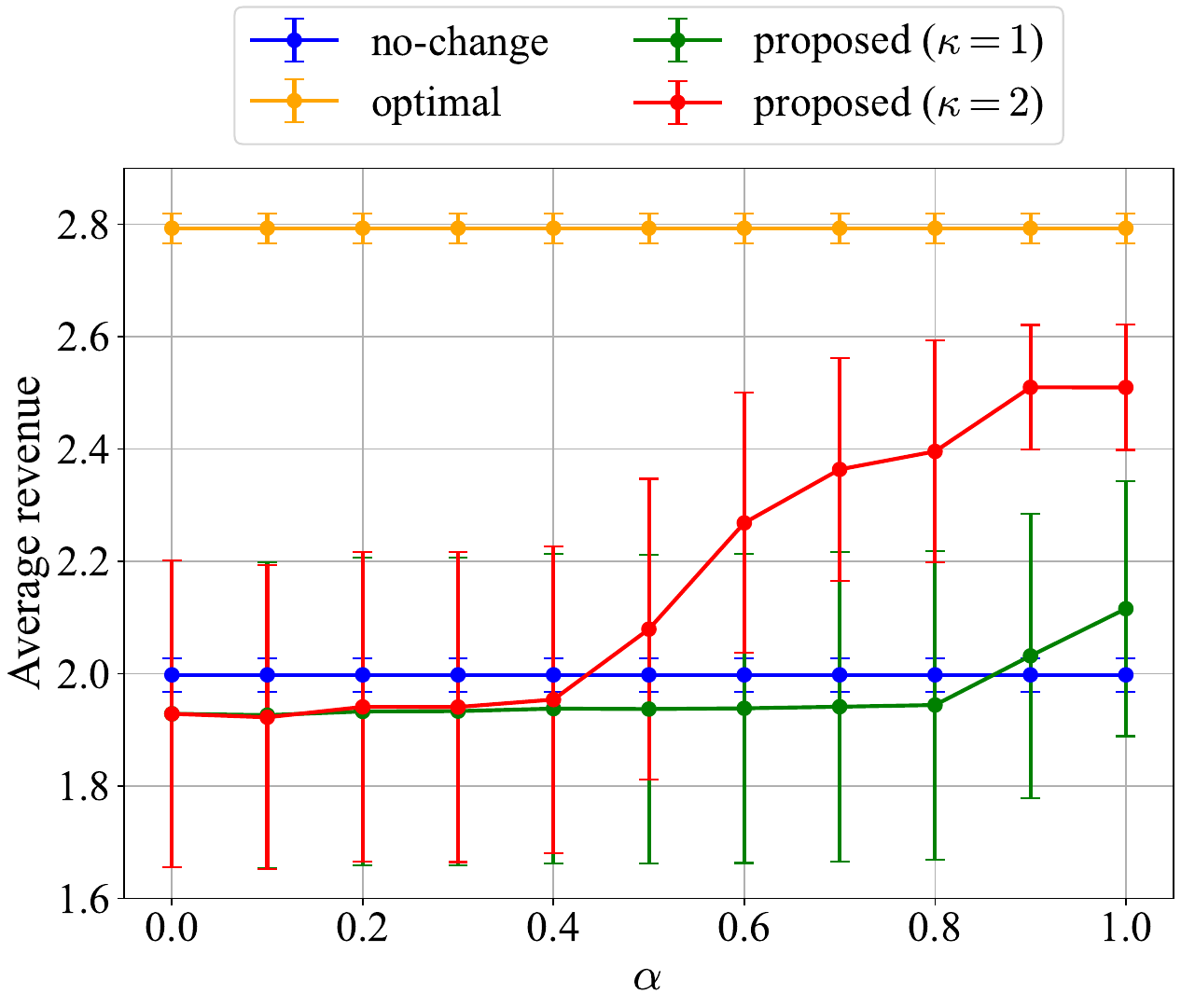}}
    \subfloat[Dataset1 ($\lvert \mathcal{I}^{{\rm train}}\rvert = 1000$)]{\includegraphics[width=2.7in, bb=0 0 607 518,clip]{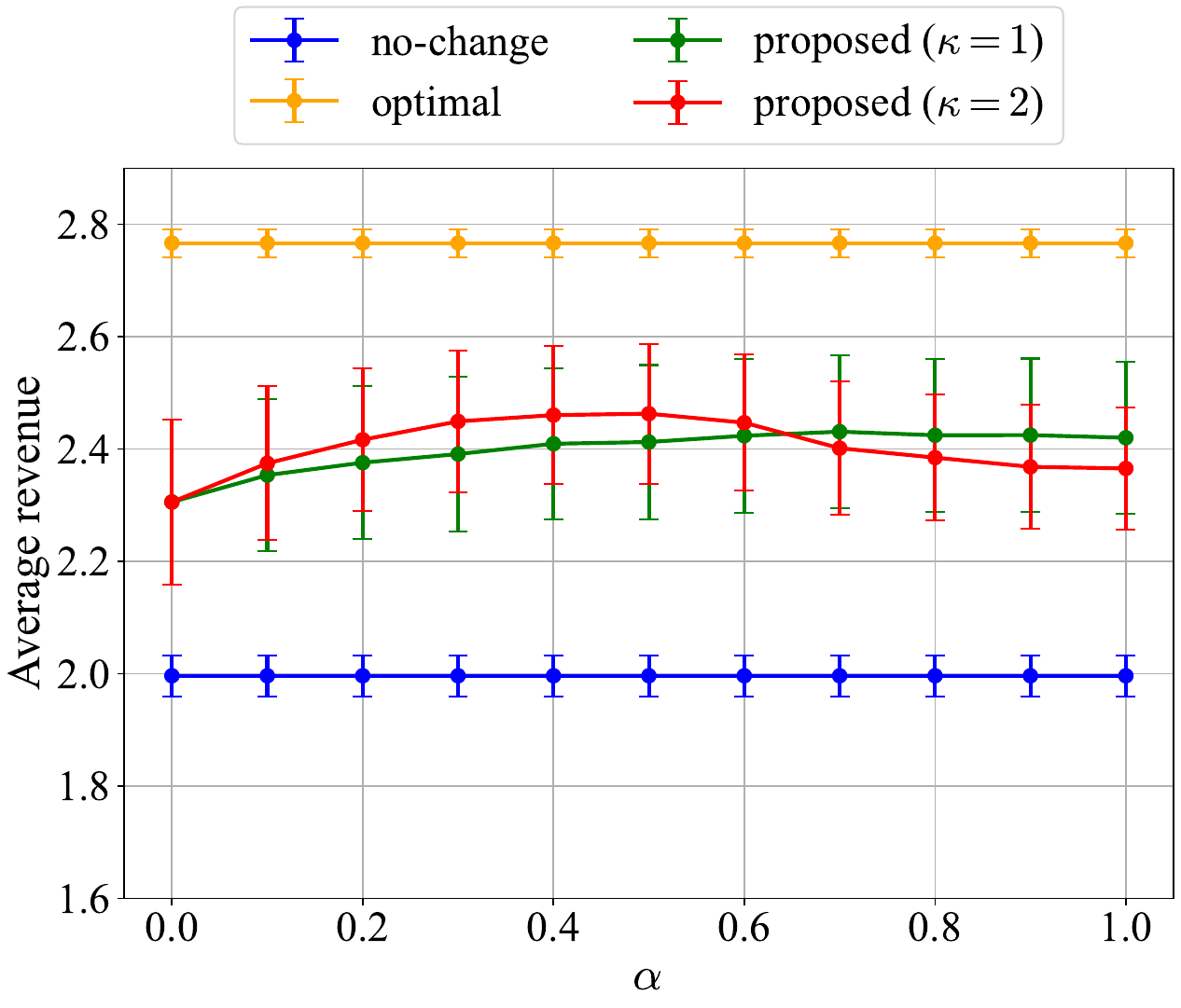}} \\
    \subfloat[Dataset2 ($\lvert \mathcal{I}^{{\rm train}}\rvert = 100$)]{\includegraphics[width=2.7in, bb=0 0 618 518,clip]{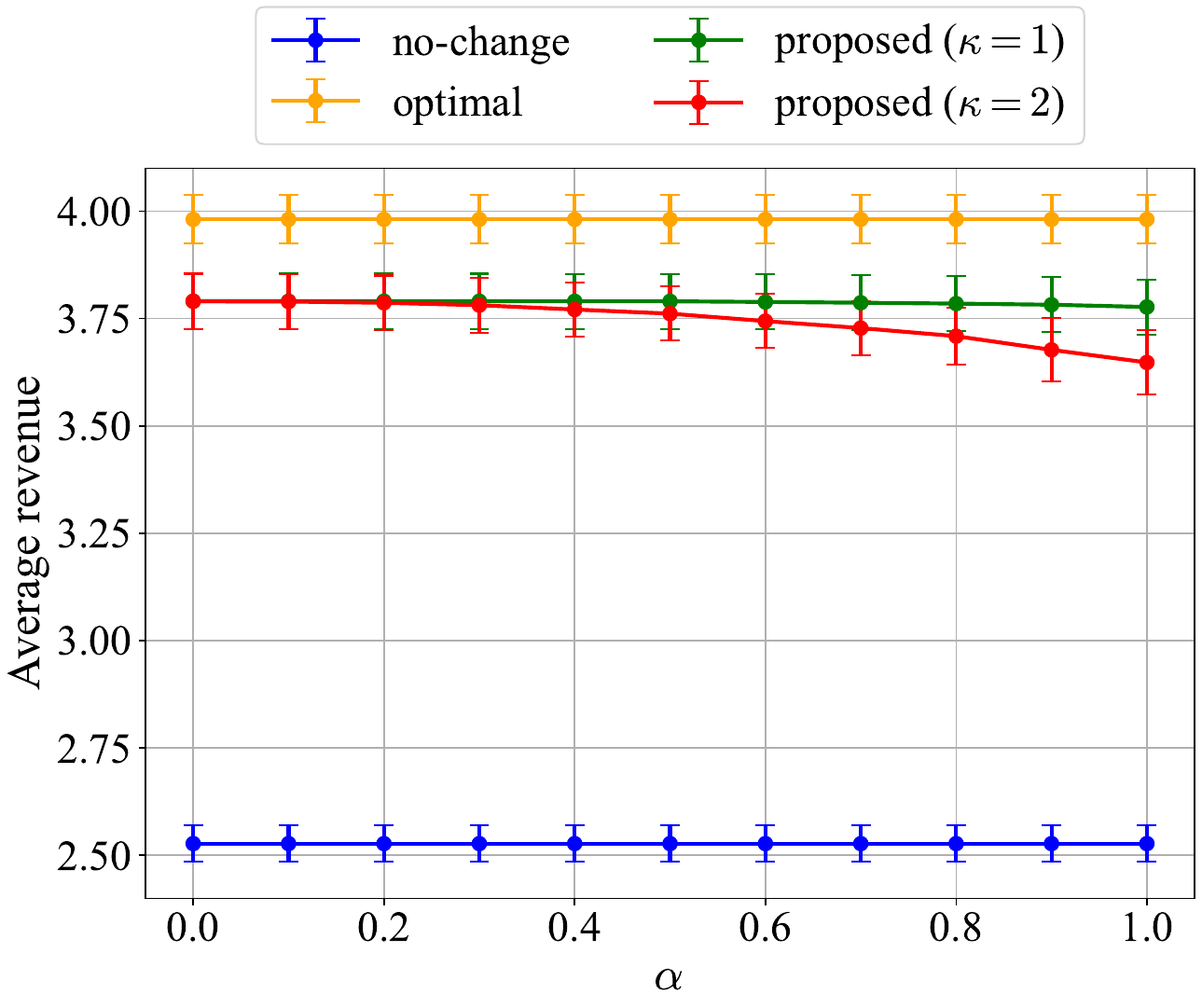}}
    \subfloat[Dataset2 ($\lvert \mathcal{I}^{{\rm train}}\rvert = 1000$)]{\includegraphics[width=2.7in, bb=0 0 618 518,clip]{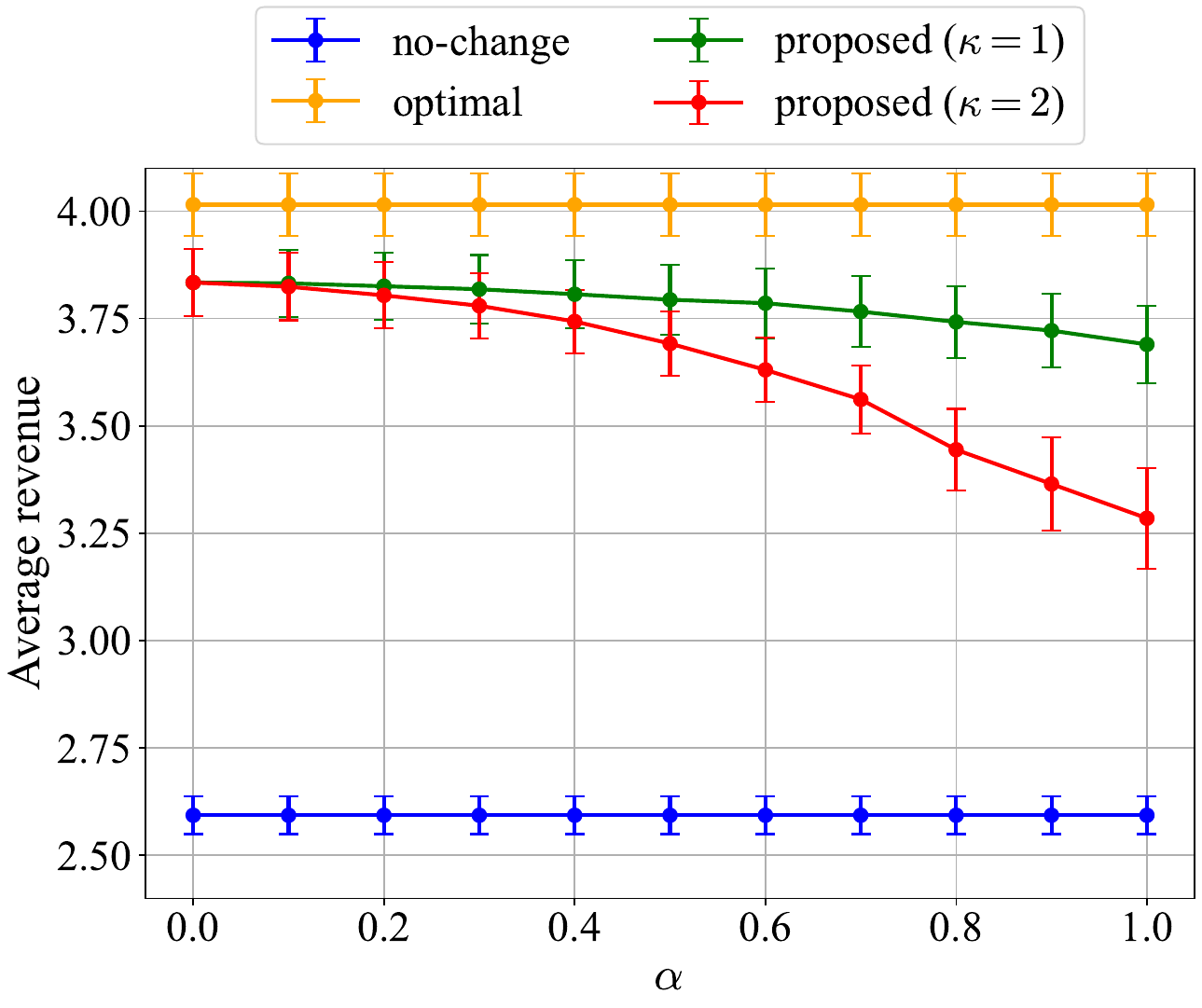}} \\
    \subfloat[Dataset3 ($\lvert \mathcal{I}^{{\rm train}}\rvert = 100$)]{\includegraphics[width=2.7in, bb=0 0 607 518,clip]{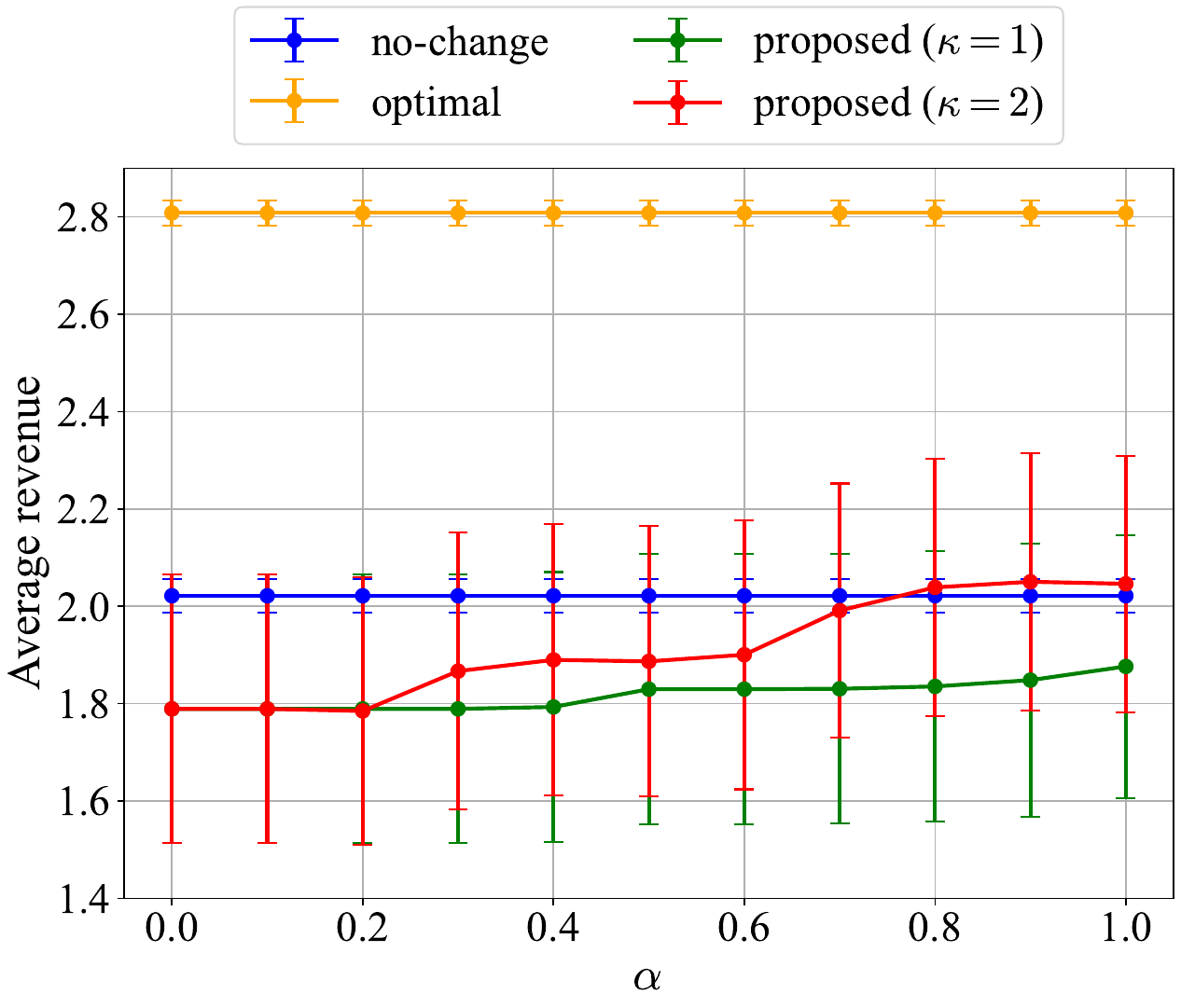}}
    \subfloat[Dataset3 ($\lvert \mathcal{I}^{{\rm train}}\rvert = 1000$)]{\includegraphics[width=2.7in, bb=0 0 607 518,clip]{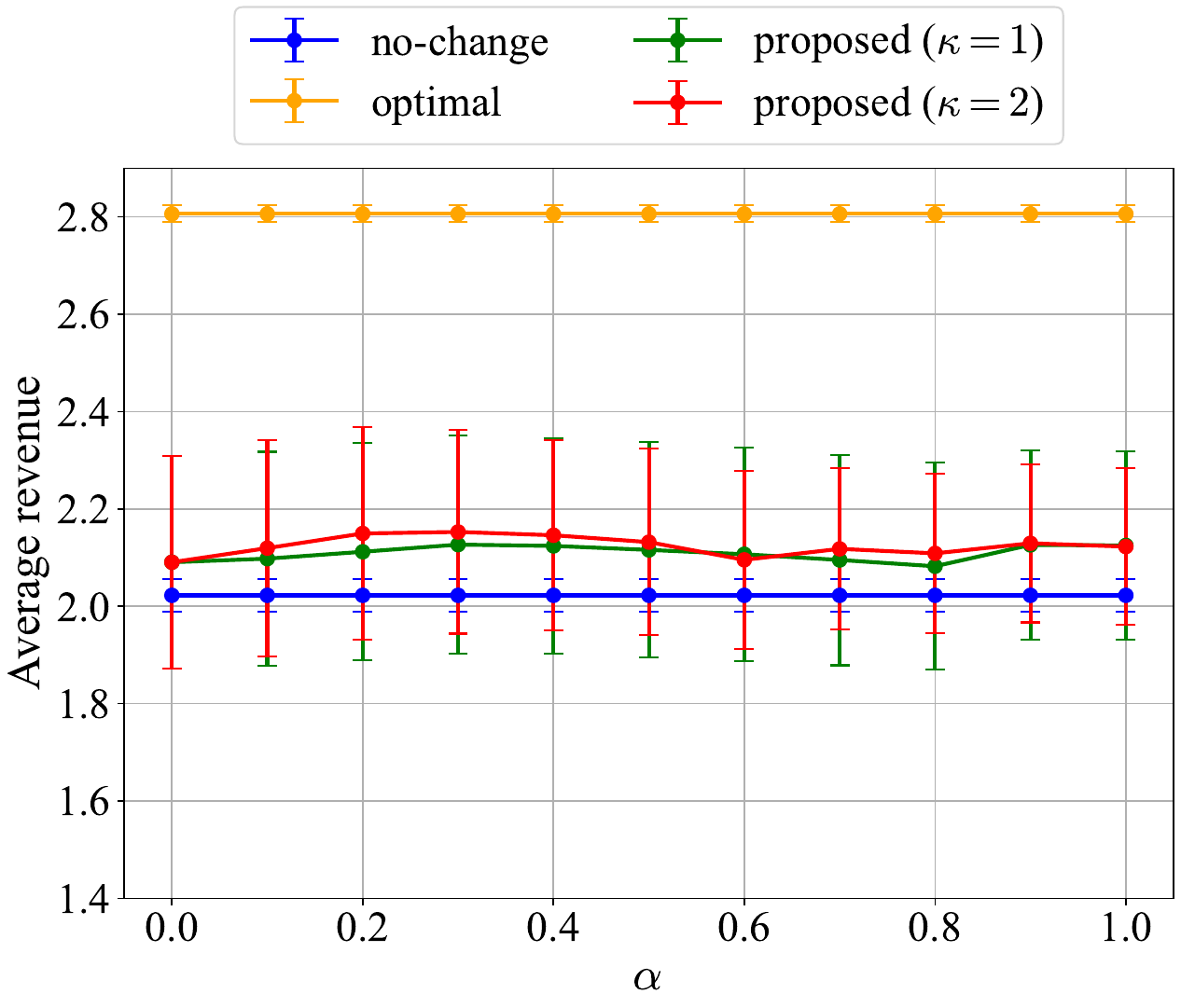}}
    \caption{Average revenues for Dataset1--3 ($\lvert \mathcal{I}^{\rm{test}} \rvert = 500$)}
    \label{fig:robust_synthetic}
\end{figure}

\begin{figure}[H]
\centering
    \subfloat[Dataset4 ($\lvert \mathcal{I}^{{\rm train}}\rvert = 100$)]{\includegraphics[width=2.7in, bb=0 0 607 518,clip]{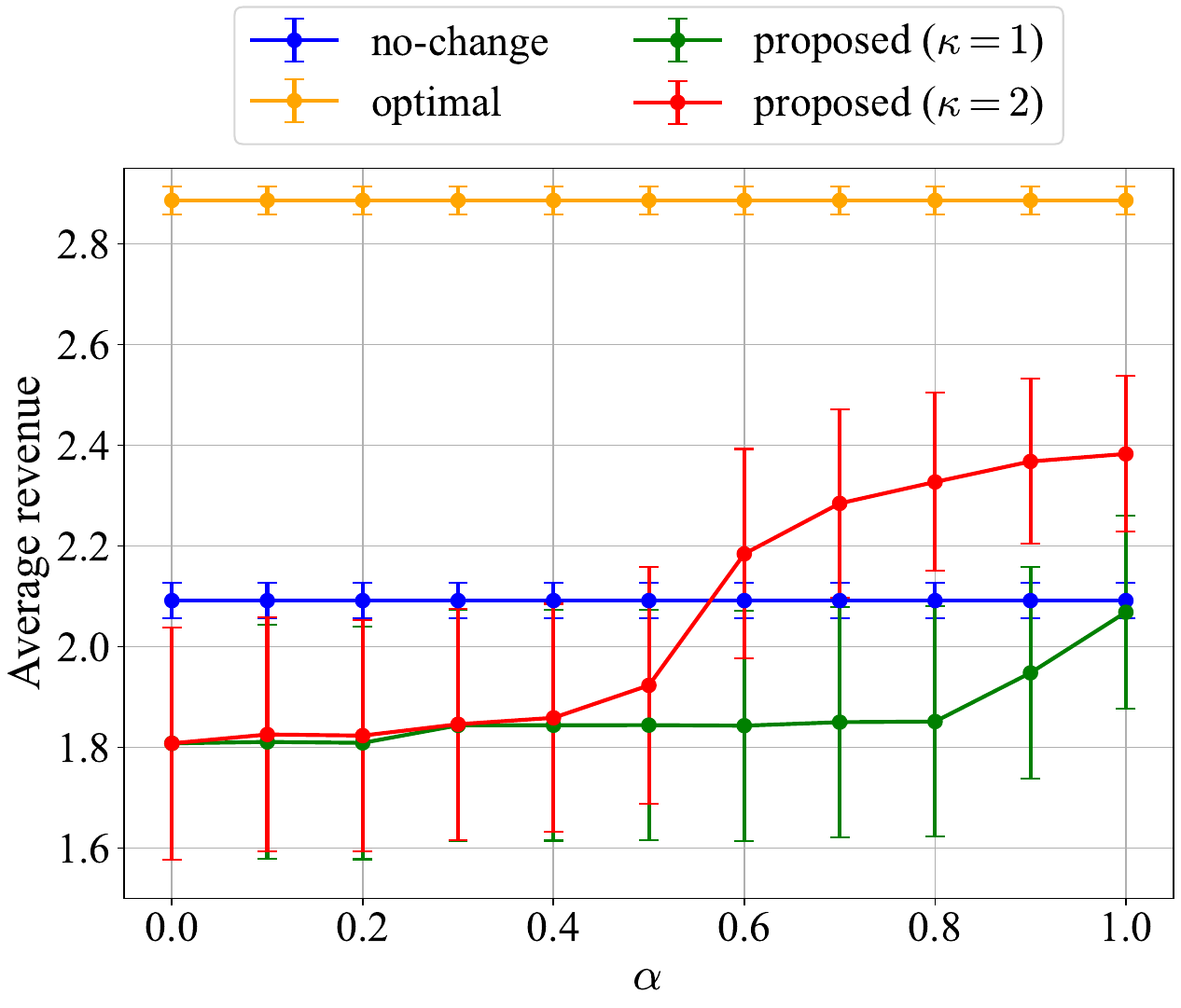}}
    \subfloat[Dataset4 ($\lvert \mathcal{I}^{{\rm train}}\rvert = 1000$)]{\includegraphics[width=2.7in, bb=0 0 607 518,clip]{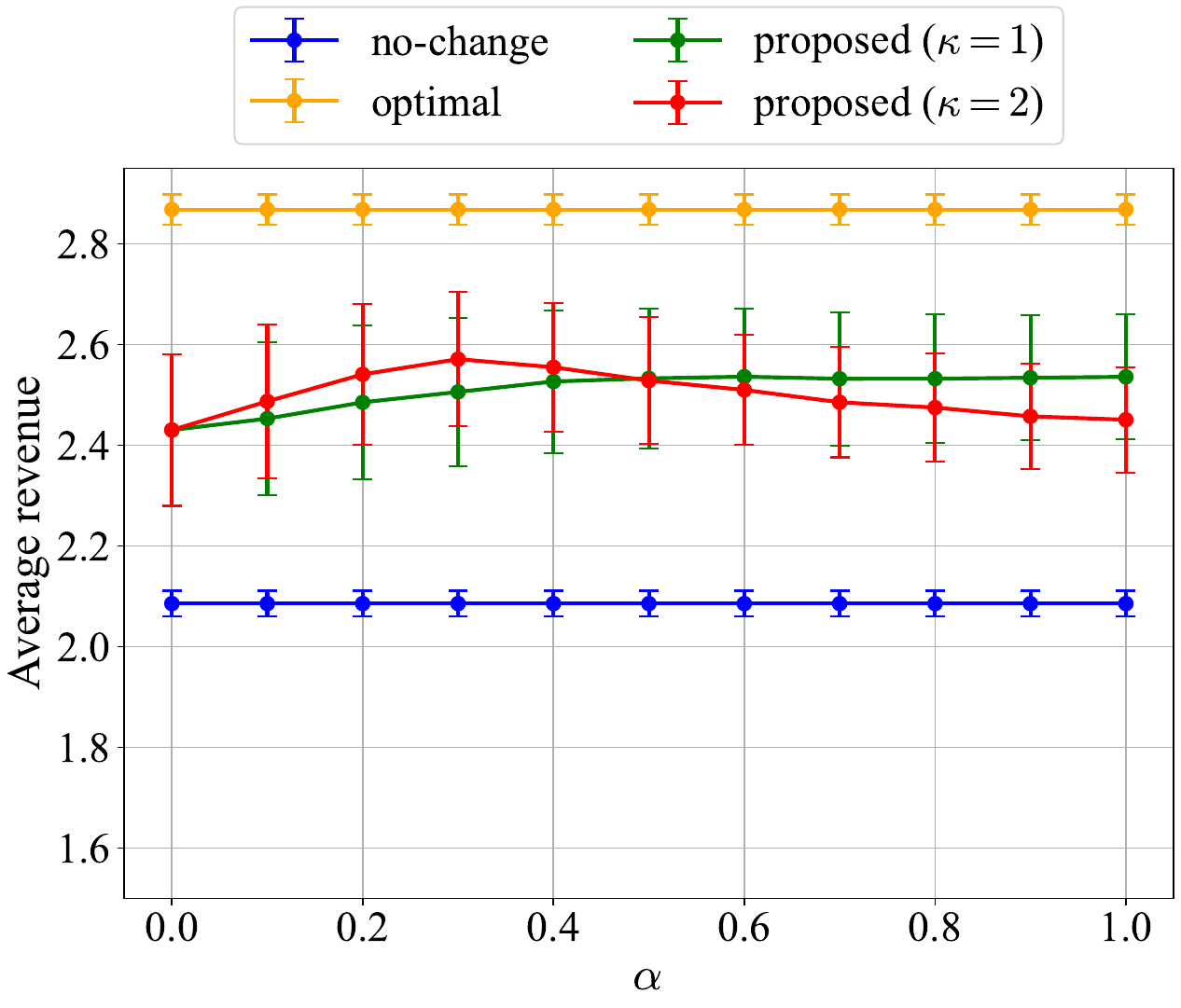}} \\
    \subfloat[Dataset5 ($\lvert \mathcal{I}^{{\rm train}}\rvert = 100$)]{\includegraphics[width=2.7in, bb=0 0 607 518,clip]{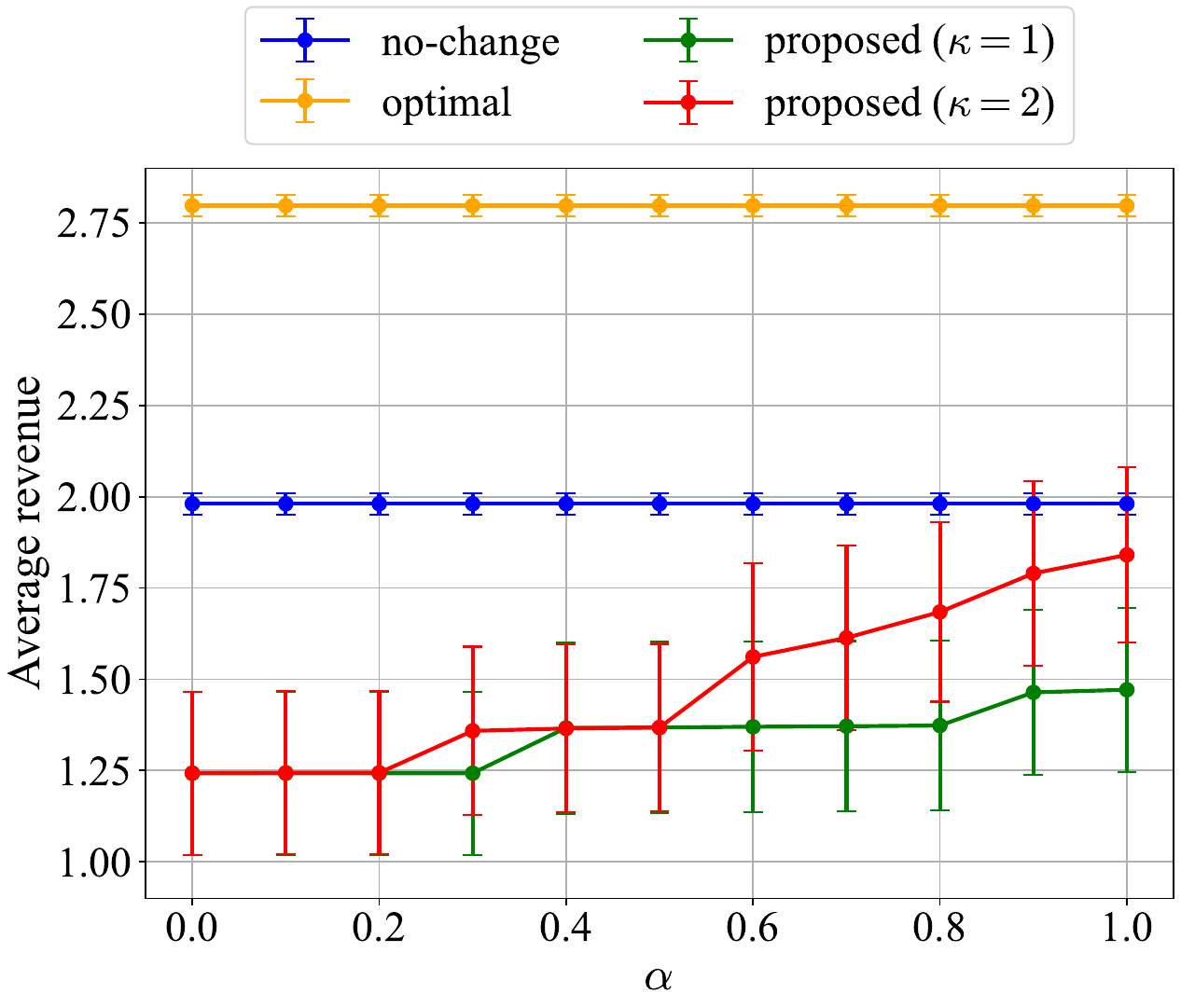}}
    \subfloat[Dataset5 ($\lvert \mathcal{I}^{{\rm train}}\rvert = 1000$)]{\includegraphics[width=2.7in, bb=0 0 607 518,clip]{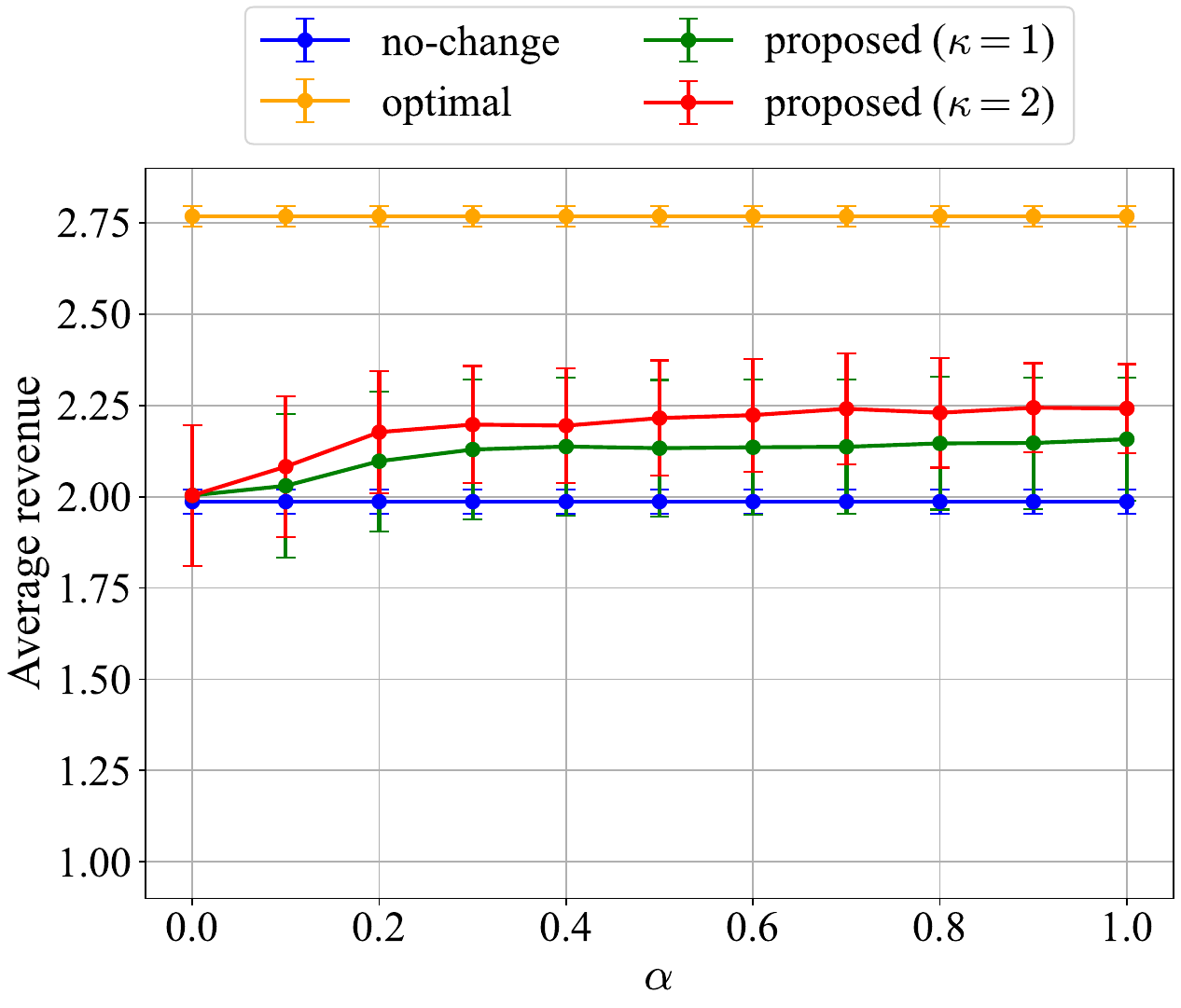}} \\
    \subfloat[Dataset6 ($\lvert \mathcal{I}^{{\rm train}}\rvert = 100$)]{\includegraphics[width=2.7in, bb=0 0 607 518,clip]{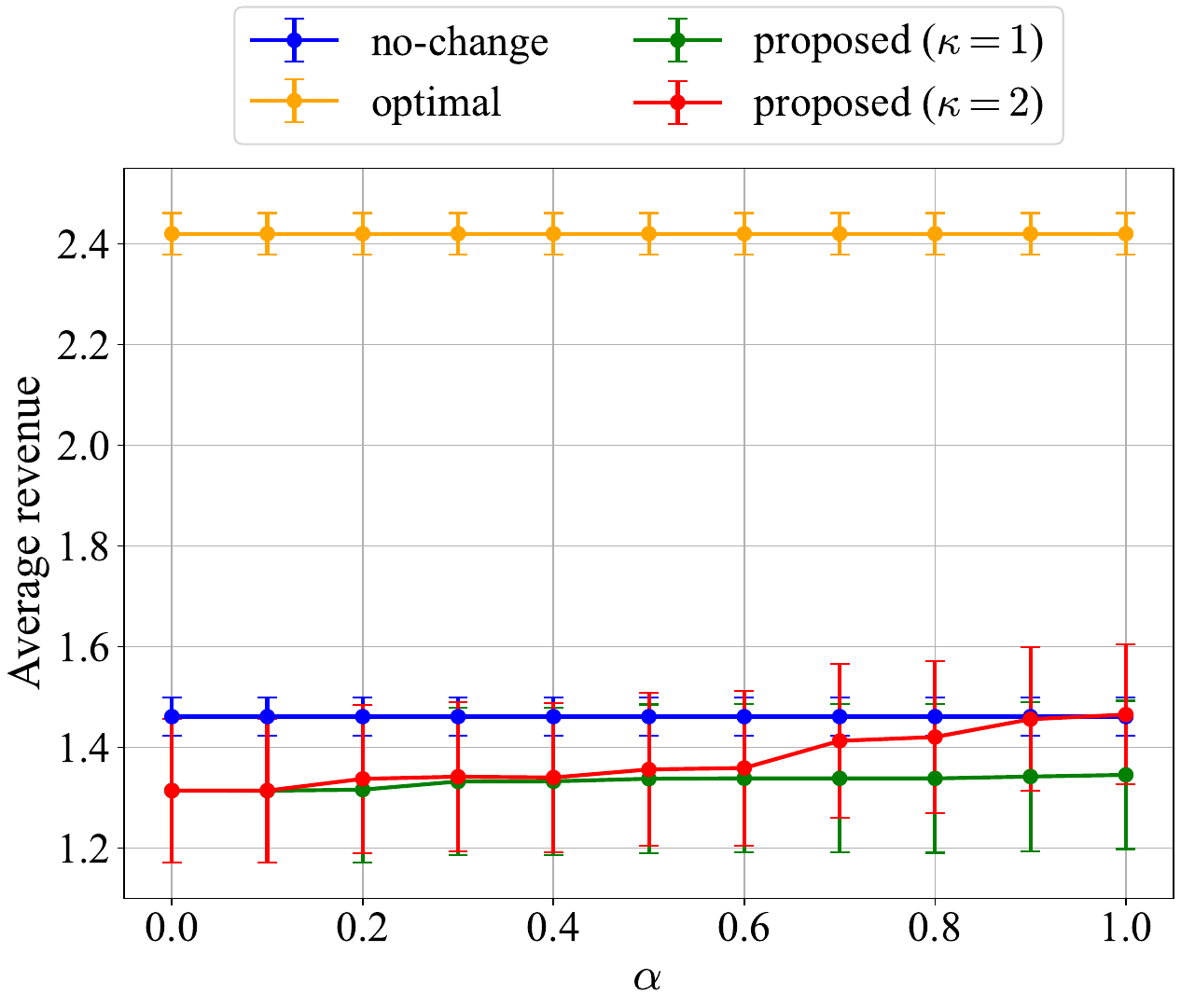}}
    \subfloat[Dataset6 ($\lvert \mathcal{I}^{{\rm train}}\rvert = 1000$)]{\includegraphics[width=2.7in, bb=0 0 607 518,clip]{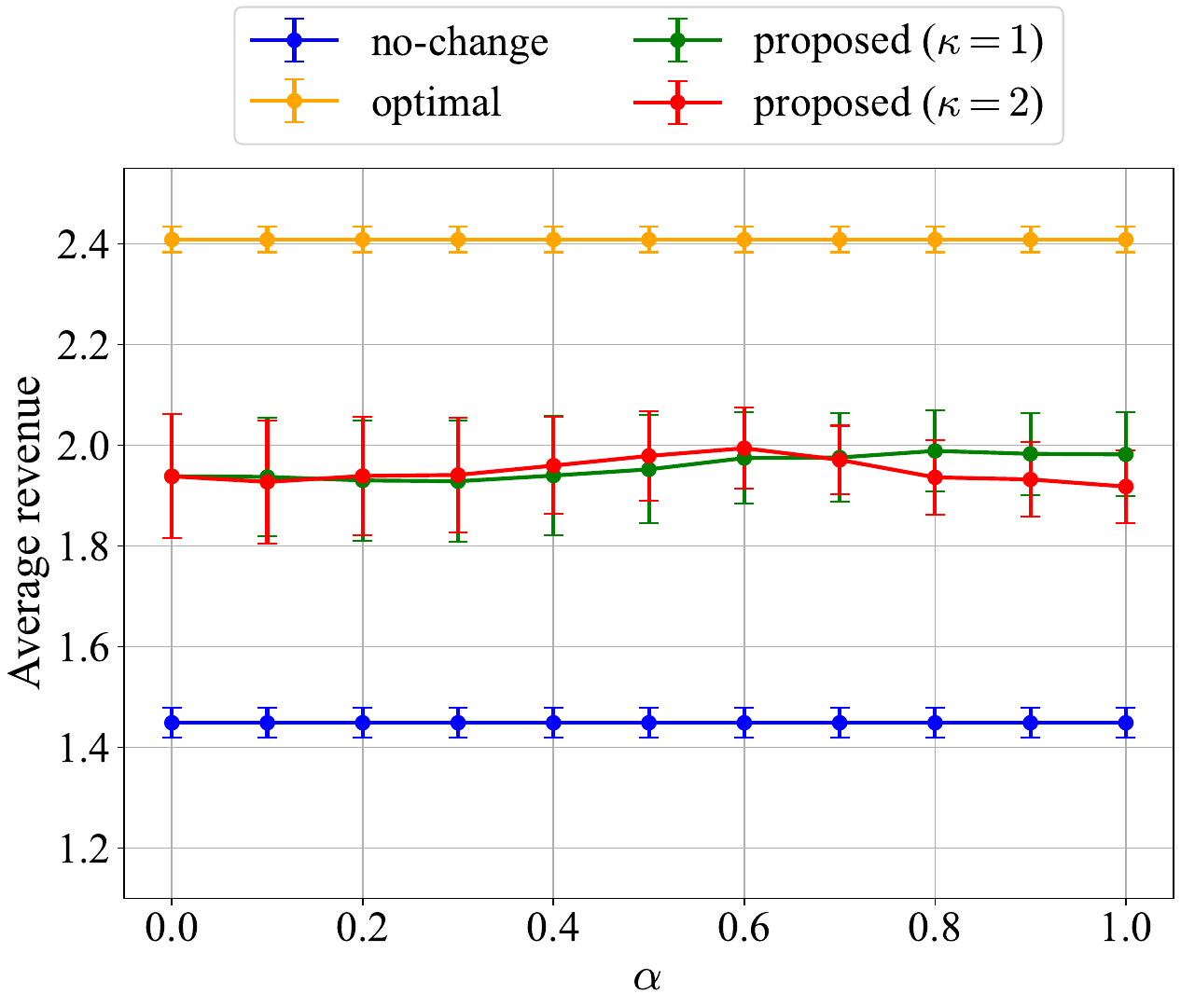}}
    \caption{Average revenues for Dataset4--6 ($\lvert \mathcal{I}^{\rm{test}} \rvert = 500$)}
    \label{fig:robust_synthetic2}
\end{figure}

Next, we examine the effect of accounting for uncertainty of predicted purchase probabilities on revenues.
Figs.~\ref{fig:robust_synthetic} and \ref{fig:robust_synthetic2} illustrate the average revenue for each dataset when varying the ratio of consumers affected by uncertainty, denoted by $\alpha$.
The legend ``optimal" represents the results for the optimal pricing under the true purchase probability model, while ``no-change" represents the results without price optimization.
The prices for the ``no-change" scenario were generated using the same purchase probability model used to produce the training data.

As shown in Figs.~\ref{fig:robust_synthetic} and \ref{fig:robust_synthetic2}, considering uncertainty improved the expected revenue in all datasets except Dataset2.
In this dataset, as shown in Figs.~\ref{fig:robust_synthetic}(c) and (d), increasing $\alpha$ led to lower revenues.
Notably, the AUC for Dataset2 remained below 0.6 despite increasing the number of consumers in the training data, indicating ineffective purchase probability prediction.
This suggests that while accounting for the prediction uncertainty may increase revenues when predictions are reasonably accurate, it may decrease revenues when prediction accuracy is poor.

In Dataset1 and Dataset4, Dataset3 and Dataset6, and Dataset5, the average revenues show different behaviors.
For Dataset1 and Dataset4, when the number of consumers in the training data was $\lvert\mathcal{I}^{{\rm train}}\rvert = 100$, lower values of $\alpha$ tended to result in average revenues below those of the no-change scenario, as illustrated in Figs.~\ref{fig:robust_synthetic}(a) and \ref{fig:robust_synthetic2}(a).
However, increasing $\alpha$ led to higher expected revenues compared to the no-change scenario, especially when the uncertainty parameter $\kappa=2$.
Conversely, with $\lvert\mathcal{I}^{{\rm train}}\rvert = 1000$, the highest expected revenues were obtained at large values of $\alpha$ for $\kappa=1$ and at smaller values of $\alpha$ for $\kappa=2$, as shown in Figs.~\ref{fig:robust_synthetic}(b) and \ref{fig:robust_synthetic2}(b).

In terms of the trends observed in Dataset3 and Dataset6, for $\lvert\mathcal{I}^{{\rm train}}\rvert = 100$, Figs.~\ref{fig:robust_synthetic}(e) and \ref{fig:robust_synthetic2}(e) demonstrate that increasing $\alpha$ tended to increase or stabilize the expected revenue, especially when $\kappa=2$.
For $\lvert\mathcal{I}^{{\rm train}}\rvert = 1000$, comparisons between Figs.~\ref{fig:robust_synthetic}(f) and \ref{fig:robust_synthetic2} (f) reveal that accounting for uncertainty (i.e., $\alpha>0$) typically resulted in slightly better expected revenues compared to scenarios that did not account for uncertainty (i.e., $\alpha=0$), although the trend was generally flat.

Regarding the trends in Dataset5, as shown in Figs.~\ref{fig:robust_synthetic2}(c) and \ref{fig:robust_synthetic2}(d), it is evident that as $\alpha$ increased, the average revenue consistently increased.

Finally, the common trends in Dataset1 and Dataset4, Dataset3 and Dataset6, and Dataset5 are as follows.
Increasing $\kappa$ could improve the expected revenue.
Furthermore, when the number of consumers in the training data was small, the effect of increasing $\alpha$ became more pronounced.
Additionally, improvements in the accuracy of predicted purchase probabilities associated with increases in $\lvert\mathcal{I}^{{\rm train}}\rvert = 1000$ contributed to higher revenues.

When the number of training data was small and the prediction accuracy was insufficient (i.e., $\lvert \mathcal{I}^{{\rm train}} \rvert = 100$), the average revenue from price optimization might decrease compared to a situation where no price optimization was performed.
Conversely, with a larger training dataset (i.e., $\lvert \mathcal{I}^{{\rm train}} \rvert = 1000$), all synthetic datasets consistently achieved higher expected revenues compared to the no-change scenario, regardless of $\alpha$.
This illustrates the critical importance of the predictive accuracy in purchase probability for effective price optimization.

Excluding Dataset2, when $\lvert \mathcal{I}^{{\rm train}}\rvert = 1000$, the area under the curve (AUC) for the testing data was above 0.8, as shown in Table \ref{tab:auc}.
This suggests that an AUC greater than 0.8 could serve as a criterion for our personalized pricing.


\subsubsection{Result for real-world dataset}
\begin{figure}[H]
\centering
\includegraphics[width=3.5in, bb=0 0 607 518,clip]{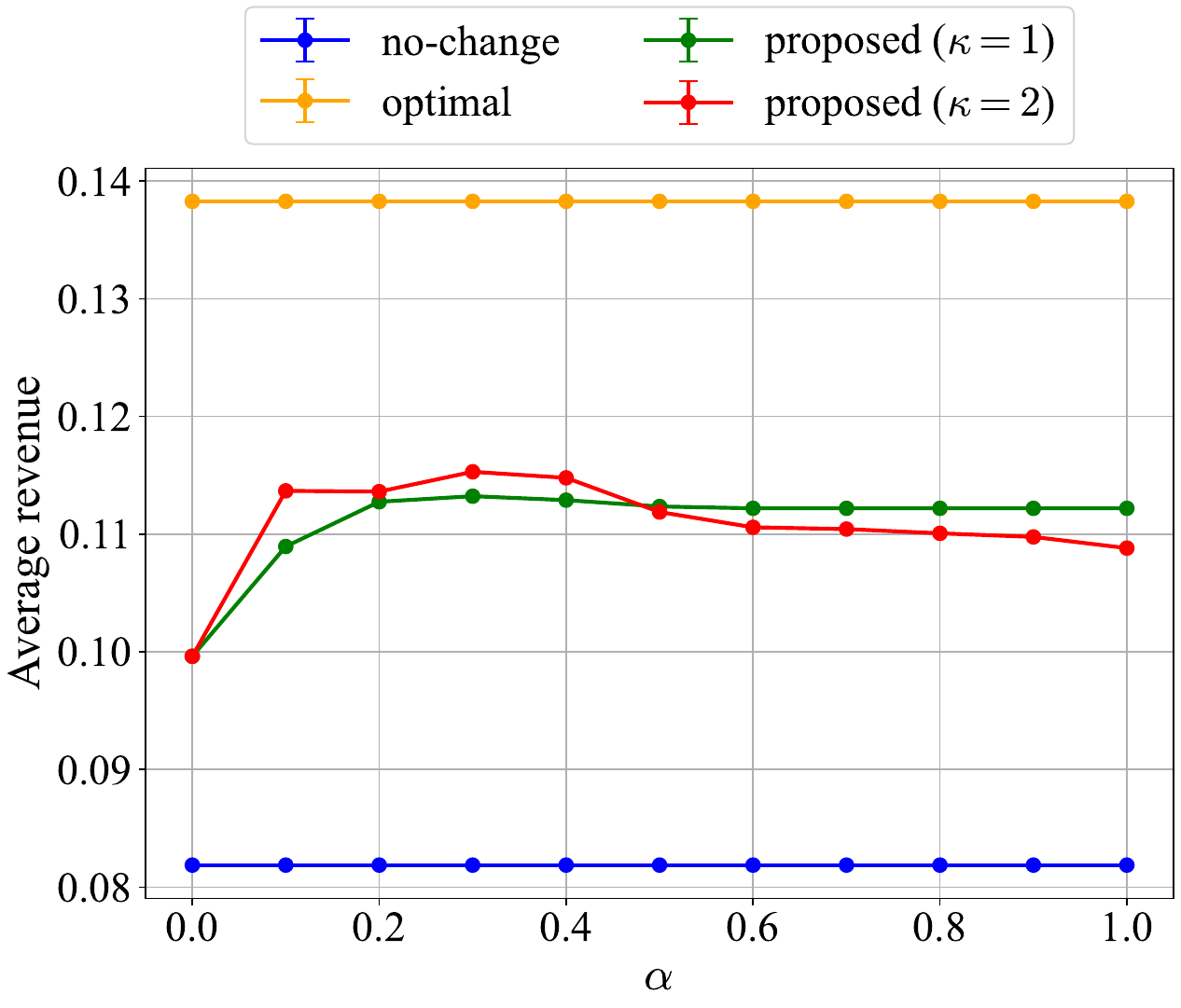}
\caption{Average revenue for ``The complete journey"}\label{robust_complete}
\end{figure}

Fig.~\ref{robust_complete} shows the average revenues for each method on ``The complete journey."
This figure demonstrates that the revenues can be improved by accounting for the uncertainty of the predicted purchase probabilities.
The accuracy of the prediction model for the purchase probabilities was ${\rm AUC} = 0.836$ for the testing data.
Regarding the computation time of the optimization problem, the exact solution could be obtained in a few seconds for all values of $\alpha$ for $\kappa=1$ and $\kappa=2$.

\subsection{Evaluation of heuristic algorithm}
Lastly, we evaluate the computational performance of our solution methods: MILO formulation~\eqref{eq:obj_robust_dual}--\eqref{eq:nu} and our decomposition algorithm (Algorithm \ref{alg2}) on synthetic datasets (Dataset1--6), since the synthetic datasets can be easily evaluated by varying to the target number of consumers $\lvert\mathcal{I}^{\rm test}\rvert$.

\subsubsection{Experimental setup}
Our heuristic algorithm allows for splitting the original problem into subproblems for any constraint $f_k({\bm z})\leq 0 ~(k \in \mathcal{K})$ on the price assignment variables ${\bm z}$.
In this experiment, we considered the following upper bound constraint on the number of price changes \citep{wang2021price,bitran2003overview}, which has been widely used in the context of price optimization:
\begin{align*}
    \sum_{i \in \mathcal{I}^{\textrm{test}}}\sum_{j \in \mathcal{J}}a_{ij}z_{ij} \leq \beta \lvert \mathcal{I}^{\textrm{test}} \rvert
\end{align*}
where $a_{ij}\in \{0, 1\}^{\lvert\mathcal{I}^{{\textrm{test}}}\times\mathcal{J}\rvert}$ is a constant that is set to 1 if price $P_{ij}$ is available to consumer $i$ and 0 otherwise.
Additionally, let $\beta \in [0, 1]$ be a parameter that specifies the number of price changes.
We assumed that the number of consumers who could set a price above the average of the price candidates was limited to less than 10\% of the total number of subjects.
Specifically, we set $a_{ij} = \{0,0,0,0,0,1,1,1,1\}_{(i, j) \in \mathcal{I}^{\rm test} \times \mathcal{J}}$ and $\beta = 0.1$.

The termination condition of our algorithm was as follows:
\[
\frac{\lvert\lvert \nabla_{{\bm \lambda}}\rho({\bm \lambda})\rvert\rvert}{t} < \varepsilon^{p}, \quad \nabla_{{\bm \lambda}}\rho({\bm \lambda}) \geq {\bm 0},
\]
where $\varepsilon^{p}$ represents the tolerance error, and we set $\varepsilon^{p}=0.01$.
The second condition was added to ensure feasibility when comparing the objective function values in feasible solutions.
The step size for updating the Lagrange multipliers was set as follows:
\[
\displaystyle{\delta^{(t)} = \frac{1}{\lvert\lvert \nabla_{{\bm \lambda}}\rho({\bm \lambda})\rvert\rvert\cdot \sqrt{t}}}.
\]
The parameter for the termination condition of the golden section search used in the proposed algorithm was set to $\varepsilon^{g} = 0.01$.
The initial value of the upper bound was set to $\max\{P_{ij}\Delta_{ij} \mid i \in \mathcal{I}^{\rm test}, j \in \mathcal{J}\}$, and the initial value of the lower bound was set to $0$.

The other parameters for the experiments were set as shown in Table \ref{tab:param}, and the time limit of the mathematical optimization solver was set to 300 seconds.
\begin{table}[H]
\centering
\caption{Parameter}
  \begin{tabular}{lcr}
    \hline
    Parameter & Value  \\
    \hline
    Number of consumers in training data $\lvert\mathcal{I}^{{\rm train}}\rvert$  & $1000$ \\
    Ratio of consumers affected by uncertainty $\alpha$  & $0.5$ \\
    Ratio of uncertainty to standard deviation $\kappa$  & $1.0$ \\
    Number of bootstrap method trials $N^{{\rm bs}}$ & $20$ \\
    \hline
  \end{tabular}
\begin{flushleft} 
\end{flushleft}
\label{tab:param}
\end{table}

\subsubsection{Computational results}
Table \ref{tab:obj_exact_heu} presents the objective function values for the exact and heuristic algorithms.
These objective function values were calculated based on the obtained solution ${\bm z}$ using the objective function of the original problem (Eq.~\eqref{eq:obj_robust_dual}).
As shown in Table \ref{tab:obj_exact_heu}, when the number of target consumers $\lvert\mathcal{I}^{\rm test}\rvert$ was small (i.e., $\lvert\mathcal{I}^{\rm test}\rvert = 100~\text{and}~\lvert\mathcal{I}^{\rm test}\rvert = 1000$), the exact algorithm yielded slightly better objective function values than the heuristic algorithm.
However, the differences were minor, and the heuristic algorithm also achieved objective function values close to those of the exact algorithm.
Furthermore, when $\lvert\mathcal{I}^{\rm test}\rvert=10000$, the heuristic algorithm showed better objective function values in Dataset1, Dataset3, Dataset4, and Dataset6.
This indicates that the revenue performance of the exact and heuristic algorithms is competitive when the number of consumers is large.

Table \ref{tab:calc_exact_heu} shows the computation time results for the exact and heuristic algorithms.
In this table, the computation time is noted as ``$>300$" if the optimal solution was not obtained within the time limit of 300 seconds in at least one of the ten trials.
The results indicate that the exact algorithm struggled to consistently find the optimal solution within the time limit when $\lvert\mathcal{I}^{\rm test}\rvert \geq 1000$.
In contrast, the heuristic algorithm completed within 300 seconds for all synthetic datasets, even when $\lvert\mathcal{I}^{\rm test}\rvert=10000$.
Recall that our algorithm can be processed in parallel, allowing for further computational speedup.
Therefore, our algorithm has proven to be fast and capable of obtaining high-quality solutions to large-scale problems.

\begin{table}[H]
    \centering
    \caption{Objective function values for exact (MILO) and heuristic (HA) algorithms}
    \scalebox{0.92}{
    \begin{tabular}{ccrrr}
        \toprule
          & & \multicolumn{3}{c}{$|\mathcal{I}^{\rm test}|$} \\
        \cmidrule(lr){3-5} 
        Dataset & Method & $100$ & $1000$ & $10000$ \\
        \midrule
        Dataset1 & MILO & $248.28~(\pm8.34)$ & $2483.18~(\pm82.51)$ & $23933.06~(\pm\phantom{0}834.31)$ \\
          & HA & $247.88~(\pm8.59)$ & $2482.98~(\pm82.53)$ & $23934.78~(\pm\phantom{0}834.46)$ \\
        \addlinespace[2mm]
        Dataset2 & MILO & $296.60~(\pm4.48)$ & $2951.66~(\pm36.42)$ & $29500.18~(\pm\phantom{0}326.58)$ \\
          & HA & $296.45~(\pm4.42)$ & $2951.70~(\pm36.42)$ & $29500.01~(\pm\phantom{0}326.77)$ \\
        \addlinespace[2mm]
        Dataset3 & MILO & $239.95~(\pm6.86)$ & $2374.76~(\pm65.50)$ & $25567.95~(\pm\phantom{0}942.01)$ \\
          & HA & $239.78~(\pm6.83)$ & $2374.55~(\pm65.48)$ & $25559.36~(\pm\phantom{0}946.15)$ \\
        \addlinespace[2mm]
        Dataset4 & MILO & $241.56~(\pm8.63)$ & $2369.35~(\pm62.22)$ & $24344.92~(\pm\phantom{0}911.17)$\\
          & HA & $241.53~(\pm8.64)$ & $2369.04~(\pm62.29)$ & $24345.61~(\pm\phantom{0}911.06)$ \\
        \addlinespace[2mm]
        Dataset5 & MILO & $236.79~(\pm6.38)$ & $2427.48~(\pm58.93)$ & $24530.26~(\pm1093.78)$\\
          & HA & $235.14~(\pm7.28)$ & $2427.20~(\pm58.97)$ & $25167.57~(\pm1009.93)$ \\
        \addlinespace[2mm]
        Dataset6 & MILO & $189.28~(\pm8.63)$ & $1714.48~(\pm65.30)$ & $18712.09~(\pm\phantom{0}711.10)$  \\
          & HA & $189.04~(\pm8.69)$ & $1714.96~(\pm65.34)$ & $18714.87~(\pm\phantom{0}710.82)$ 
        \\ \bottomrule
    \end{tabular}
    }
    \label{tab:obj_exact_heu}
\end{table}

\begin{table}[H]
    \centering
        \caption{Computation time (s) for exact (MILO) and heuristic (HA) algorithms}
    \catcode`?=\active \def?{\phantom{0}}% ?:=\phantom{0}
    \begin{tabular}{ccrrr}
        \toprule
          & & \multicolumn{3}{c}{$|\mathcal{I}^{\rm test}|$} \\
        \cmidrule(lr){3-5} 
        Dataset & Method & $100$ & $1000$ & $10000$ \\
        \midrule
        Dataset1 & MILO & $0.24~(\pm0.09)$ & $>300$ & $>300$ \\
          & HA & $0.16~(\pm0.07)$ & $2.35~(\pm?0.81)$ & $13.25~(\pm3.97)$ \\
        \addlinespace[2mm]
        Dataset2 & MILO & $0.09~(\pm0.02)$ & $14.81~(\pm?4.46)$ & $>300$ \\
          & HA & $0.01~(\pm0.00)$ & $0.18~(\pm?0.04)$ & $1.79~(\pm0.71)$ \\
        \addlinespace[2mm]
        Dataset3 & MILO & $0.33~(\pm0.10)$ & $>300$ & $>300$ \\
          & HA & $0.07~(\pm0.02)$ & $1.58~(\pm?0.87)$ & $9.88~(\pm4.22)$ \\
        \addlinespace[2mm]
        Dataset4 & MILO & $0.32~(\pm0.13)$ & $>300$ & $>300$\\
          & HA & $0.03~(\pm0.01)$ & $1.69~(\pm?0.81)$ & $6.72~(\pm3.18)$ \\
        \addlinespace[2mm]
        Dataset5 & MILO & $6.32~(\pm5.98)$ & $40.72~(\pm12.79)$ & $>300$\\
          & HA & $0.08~(\pm0.05)$ & $3.23~(\pm?1.40)$ & $15.70~(\pm5.52)$ \\
        \addlinespace[2mm]
        Dataset6 & MILO & $0.14~(\pm0.04)$ & $>300$ & $>300$  \\
          & HA & $0.07~(\pm0.04)$ & $3.23~(\pm?0.88)$ & $22.48~(\pm6.31)$ 
    \\ \bottomrule
    \end{tabular}
\label{tab:calc_exact_heu}
\end{table}

\section{Conclusion}
We proposed a robust optimization model for personalized pricing that accounts for the uncertainty of predicted purchase probabilities for a single item.
Specifically, we calculated the uncertainty of predicted purchase probabilities using the bootstrap method and formulated the optimization problem as a mixed-integer linear optimization problem, which can be solved exactly using mathematical optimization solvers.
Moreover, we developed a scalable heuristic algorithm to efficiently obtain high-quality solutions for large-scale problems.

We conducted numerical experiments using both synthetic and real-world datasets and confirmed the improvement in expected revenues by considering the uncertainty of the predicted purchase probabilities.
Furthermore, we demonstrated our heuristic algorithm can quickly obtain high-quality solutions for large-scale problems.

Future directions for this study include extending the pricing model to multiple items by considering the estimation of cross-elasticity among items \citep{ito2017optimization, ikeda2023prescriptive} and applying appropriate price ranges \citep{ikeda2023operating, ikeda2024interpretable} for each consumer.
A further possible extension is to develop a method for dealing with uncertainty, similar to portfolio optimization approaches \citep{yanagi2024robust, uehara2024robust}.
% Additionally, verifying the performance of this study in actual service operations remains a future challenge.




\bibliographystyle{elsarticle-harv}
% \bibliographystyle{cas-model2-names}

\bibliography{ref}

\begin{thebibliography}{37}
\expandafter\ifx\csname natexlab\endcsname\relax\def\natexlab#1{#1}\fi
\providecommand{\url}[1]{\texttt{#1}}
\providecommand{\href}[2]{#2}
\providecommand{\path}[1]{#1}
\providecommand{\DOIprefix}{doi:}
\providecommand{\ArXivprefix}{arXiv:}
\providecommand{\URLprefix}{URL: }
\providecommand{\Pubmedprefix}{pmid:}
\providecommand{\doi}[1]{\href{http://dx.doi.org/#1}{\path{#1}}}
\providecommand{\Pubmed}[1]{\href{pmid:#1}{\path{#1}}}
\providecommand{\bibinfo}[2]{#2}
\ifx\xfnm\relax \def\xfnm[#1]{\unskip,\space#1}\fi
%Type = Article
\bibitem[{Allenby and Rossi(1998)}]{allenby1998marketing}
\bibinfo{author}{Allenby, G.M.}, \bibinfo{author}{Rossi, P.E.},
  \bibinfo{year}{1998}.
\newblock \bibinfo{title}{Marketing models of consumer heterogeneity}.
\newblock \bibinfo{journal}{Journal of econometrics} \bibinfo{volume}{89},
  \bibinfo{pages}{57--78}.
%Type = Article
\bibitem[{Amram et~al.(2022)Amram, Dunn and Zhuo}]{amram2022optimal}
\bibinfo{author}{Amram, M.}, \bibinfo{author}{Dunn, J.}, \bibinfo{author}{Zhuo,
  Y.D.}, \bibinfo{year}{2022}.
\newblock \bibinfo{title}{Optimal policy trees}.
\newblock \bibinfo{journal}{Machine Learning} \bibinfo{volume}{111},
  \bibinfo{pages}{2741--2768}.
%Type = Article
\bibitem[{Athey and Imbens(2016)}]{athey2016recursive}
\bibinfo{author}{Athey, S.}, \bibinfo{author}{Imbens, G.},
  \bibinfo{year}{2016}.
\newblock \bibinfo{title}{Recursive partitioning for heterogeneous causal
  effects}.
\newblock \bibinfo{journal}{Proceedings of the National Academy of Sciences}
  \bibinfo{volume}{113}, \bibinfo{pages}{7353--7360}.
%Type = Article
\bibitem[{Bertsimas and Sim(2004)}]{bertsimas2004price}
\bibinfo{author}{Bertsimas, D.}, \bibinfo{author}{Sim, M.},
  \bibinfo{year}{2004}.
\newblock \bibinfo{title}{The price of robustness}.
\newblock \bibinfo{journal}{Operations Research} \bibinfo{volume}{52},
  \bibinfo{pages}{35--53}.
%Type = Article
\bibitem[{Bertsimas and Van~Parys(2022)}]{bertsimas2022bootstrap}
\bibinfo{author}{Bertsimas, D.}, \bibinfo{author}{Van~Parys, B.},
  \bibinfo{year}{2022}.
\newblock \bibinfo{title}{Bootstrap robust prescriptive analytics}.
\newblock \bibinfo{journal}{Mathematical Programming} \bibinfo{volume}{195},
  \bibinfo{pages}{39--78}.
%Type = Inproceedings
\bibitem[{Biggs et~al.(2021)Biggs, Sun and Ettl}]{biggs2021model}
\bibinfo{author}{Biggs, M.}, \bibinfo{author}{Sun, W.}, \bibinfo{author}{Ettl,
  M.}, \bibinfo{year}{2021}.
\newblock \bibinfo{title}{Model distillation for revenue optimization:
  Interpretable personalized pricing}, in: \bibinfo{booktitle}{International
  Conference on Machine Learning}, \bibinfo{organization}{PMLR}. pp.
  \bibinfo{pages}{946--956}.
%Type = Article
\bibitem[{Bitran and Caldentey(2003)}]{bitran2003overview}
\bibinfo{author}{Bitran, G.}, \bibinfo{author}{Caldentey, R.},
  \bibinfo{year}{2003}.
\newblock \bibinfo{title}{An overview of pricing models for revenue
  management}.
\newblock \bibinfo{journal}{Manufacturing \& Service Operations Management}
  \bibinfo{volume}{5}, \bibinfo{pages}{203--229}.
%Type = Article
\bibitem[{Boyd et~al.(2003)Boyd, Xiao and Mutapcic}]{boyd2003subgradient}
\bibinfo{author}{Boyd, S.}, \bibinfo{author}{Xiao, L.},
  \bibinfo{author}{Mutapcic, A.}, \bibinfo{year}{2003}.
\newblock \bibinfo{title}{Subgradient methods}.
\newblock \bibinfo{journal}{Lecture Notes of EE392o, Stanford University,
  Autumn Quarter} \bibinfo{volume}{2004}.
%Type = Article
\bibitem[{Chen et~al.(2023a)Chen, Cire, Hu and Lagzi}]{chen2023model}
\bibinfo{author}{Chen, N.}, \bibinfo{author}{Cire, A.A.}, \bibinfo{author}{Hu,
  M.}, \bibinfo{author}{Lagzi, S.}, \bibinfo{year}{2023}a.
\newblock \bibinfo{title}{Model-free assortment pricing with transaction data}.
\newblock \bibinfo{journal}{Management Science} \bibinfo{volume}{69},
  \bibinfo{pages}{5830--5847}.
%Type = Inproceedings
\bibitem[{Chen et~al.(2023b)Chen, Xu, Zhao and Zhou}]{chen2023personalized}
\bibinfo{author}{Chen, X.}, \bibinfo{author}{Xu, Z.}, \bibinfo{author}{Zhao,
  Z.}, \bibinfo{author}{Zhou, Y.}, \bibinfo{year}{2023}b.
\newblock \bibinfo{title}{Personalized pricing with group fairness constraint},
  in: \bibinfo{booktitle}{Proceedings of the 2023 ACM Conference on Fairness,
  Accountability, and Transparency}, pp. \bibinfo{pages}{1520--1530}.
%Type = Article
\bibitem[{Delage and Ye(2010)}]{delage2010distributionally}
\bibinfo{author}{Delage, E.}, \bibinfo{author}{Ye, Y.}, \bibinfo{year}{2010}.
\newblock \bibinfo{title}{Distributionally robust optimization under moment
  uncertainty with application to data-driven problems}.
\newblock \bibinfo{journal}{Operations research} \bibinfo{volume}{58},
  \bibinfo{pages}{595--612}.
%Type = Article
\bibitem[{Den~Boer and Zwart(2014)}]{den2014simultaneously}
\bibinfo{author}{Den~Boer, A.V.}, \bibinfo{author}{Zwart, B.},
  \bibinfo{year}{2014}.
\newblock \bibinfo{title}{Simultaneously learning and optimizing using
  controlled variance pricing}.
\newblock \bibinfo{journal}{Management science} \bibinfo{volume}{60},
  \bibinfo{pages}{770--783}.
%Type = Article
\bibitem[{Duchi et~al.(2021)Duchi, Glynn and Namkoong}]{duchi2021statistics}
\bibinfo{author}{Duchi, J.C.}, \bibinfo{author}{Glynn, P.W.},
  \bibinfo{author}{Namkoong, H.}, \bibinfo{year}{2021}.
\newblock \bibinfo{title}{Statistics of robust optimization: A generalized
  empirical likelihood approach}.
\newblock \bibinfo{journal}{Mathematics of Operations Research}
  \bibinfo{volume}{46}, \bibinfo{pages}{946--969}.
%Type = Incollection
\bibitem[{Efron(1992)}]{efron1992bootstrap}
\bibinfo{author}{Efron, B.}, \bibinfo{year}{1992}.
\newblock \bibinfo{title}{Bootstrap methods: Another look at the jackknife},
  in: \bibinfo{booktitle}{Breakthroughs in Statistics: Methodology and
  Distribution}. \bibinfo{publisher}{Springer}, pp. \bibinfo{pages}{569--593}.
%Type = Article
\bibitem[{Elmachtoub et~al.(2021)Elmachtoub, Gupta and
  Hamilton}]{elmachtoub2021value}
\bibinfo{author}{Elmachtoub, A.N.}, \bibinfo{author}{Gupta, V.},
  \bibinfo{author}{Hamilton, M.L.}, \bibinfo{year}{2021}.
\newblock \bibinfo{title}{The value of personalized pricing}.
\newblock \bibinfo{journal}{Management Science} \bibinfo{volume}{67},
  \bibinfo{pages}{6055--6070}.
%Type = Article
\bibitem[{Elreedy et~al.(2021)Elreedy, Atiya and Shaheen}]{elreedy2021novel}
\bibinfo{author}{Elreedy, D.}, \bibinfo{author}{Atiya, A.F.},
  \bibinfo{author}{Shaheen, S.I.}, \bibinfo{year}{2021}.
\newblock \bibinfo{title}{Novel pricing strategies for revenue maximization and
  demand learning using an exploration--exploitation framework}.
\newblock \bibinfo{journal}{Soft Computing} \bibinfo{volume}{25},
  \bibinfo{pages}{11711--11733}.
%Type = Article
\bibitem[{Feldman et~al.(2015)Feldman, Trzcinka and Winer}]{feldman2015pricing}
\bibinfo{author}{Feldman, D.}, \bibinfo{author}{Trzcinka, C.},
  \bibinfo{author}{Winer, R.S.}, \bibinfo{year}{2015}.
\newblock \bibinfo{title}{Pricing under noisy signaling}.
\newblock \bibinfo{journal}{Review of Quantitative Finance and Accounting}
  \bibinfo{volume}{45}, \bibinfo{pages}{435--454}.
%Type = Incollection
\bibitem[{Geoffrion(2009)}]{geoffrion2009lagrangean}
\bibinfo{author}{Geoffrion, A.M.}, \bibinfo{year}{2009}.
\newblock \bibinfo{title}{Lagrangean relaxation for integer programming}, in:
  \bibinfo{booktitle}{Approaches to integer programming}.
  \bibinfo{publisher}{Springer}, pp. \bibinfo{pages}{82--114}.
%Type = Article
\bibitem[{Golrezaei et~al.(2014)Golrezaei, Nazerzadeh and
  Rusmevichientong}]{golrezaei2014real}
\bibinfo{author}{Golrezaei, N.}, \bibinfo{author}{Nazerzadeh, H.},
  \bibinfo{author}{Rusmevichientong, P.}, \bibinfo{year}{2014}.
\newblock \bibinfo{title}{Real-time optimization of personalized assortments}.
\newblock \bibinfo{journal}{Management Science} \bibinfo{volume}{60},
  \bibinfo{pages}{1532--1551}.
%Type = Article
\bibitem[{Guignard(2003)}]{guignard2003lagrangean}
\bibinfo{author}{Guignard, M.}, \bibinfo{year}{2003}.
\newblock \bibinfo{title}{Lagrangean relaxation}.
\newblock \bibinfo{journal}{Top} \bibinfo{volume}{11},
  \bibinfo{pages}{151--200}.
%Type = Article
\bibitem[{Ikeda et~al.(2023a)Ikeda, Nishimura, Sukegawa and
  Takano}]{ikeda2023prescriptive}
\bibinfo{author}{Ikeda, S.}, \bibinfo{author}{Nishimura, N.},
  \bibinfo{author}{Sukegawa, N.}, \bibinfo{author}{Takano, Y.},
  \bibinfo{year}{2023}a.
\newblock \bibinfo{title}{Prescriptive price optimization using optimal
  regression trees}.
\newblock \bibinfo{journal}{Operations Research Perspectives}
  \bibinfo{volume}{11}, \bibinfo{pages}{100290}.
%Type = Inproceedings
\bibitem[{Ikeda et~al.(2023b)Ikeda, Nishimura and Umetani}]{ikeda2023operating}
\bibinfo{author}{Ikeda, S.}, \bibinfo{author}{Nishimura, N.},
  \bibinfo{author}{Umetani, S.}, \bibinfo{year}{2023}b.
\newblock \bibinfo{title}{Operating range estimation for price optimization},
  in: \bibinfo{booktitle}{2023 IEEE International Conference on Data Mining
  Workshops (ICDMW)}, \bibinfo{organization}{IEEE}. pp.
  \bibinfo{pages}{31--36}.
%Type = Article
\bibitem[{Ikeda et~al.(2024)Ikeda, Nishimura and
  Umetani}]{ikeda2024interpretable}
\bibinfo{author}{Ikeda, S.}, \bibinfo{author}{Nishimura, N.},
  \bibinfo{author}{Umetani, S.}, \bibinfo{year}{2024}.
\newblock \bibinfo{title}{Interpretable price bounds estimation with shape
  constraints in price optimization}.
\newblock \bibinfo{journal}{arXiv preprint arXiv:2405.14909} .
%Type = Inproceedings
\bibitem[{Ito and Fujimaki(2017)}]{ito2017optimization}
\bibinfo{author}{Ito, S.}, \bibinfo{author}{Fujimaki, R.},
  \bibinfo{year}{2017}.
\newblock \bibinfo{title}{Optimization beyond prediction: Prescriptive price
  optimization}, in: \bibinfo{booktitle}{Proceedings of the 23rd ACM SIGKDD
  international conference on knowledge discovery and data mining}, pp.
  \bibinfo{pages}{1833--1841}.
%Type = Article
\bibitem[{Jin et~al.(2021)Jin, Lin and Zhou}]{jin2021price}
\bibinfo{author}{Jin, Q.}, \bibinfo{author}{Lin, J.Y.}, \bibinfo{author}{Zhou,
  S.X.}, \bibinfo{year}{2021}.
\newblock \bibinfo{title}{Price discounts and personalized product assortments
  under multinomial logit choice model: A robust approach}.
\newblock \bibinfo{journal}{IISE Transactions} \bibinfo{volume}{53},
  \bibinfo{pages}{453--471}.
%Type = Inproceedings
\bibitem[{Kallus(2017)}]{kallus2017recursive}
\bibinfo{author}{Kallus, N.}, \bibinfo{year}{2017}.
\newblock \bibinfo{title}{Recursive partitioning for personalization using
  observational data}, in: \bibinfo{booktitle}{International Conference on
  Machine Learning}, \bibinfo{organization}{PMLR}. pp.
  \bibinfo{pages}{1789--1798}.
%Type = Inproceedings
\bibitem[{Kallus and Zhou(2021)}]{kallus2021fairness}
\bibinfo{author}{Kallus, N.}, \bibinfo{author}{Zhou, A.}, \bibinfo{year}{2021}.
\newblock \bibinfo{title}{Fairness, welfare, and equity in personalized
  pricing}, in: \bibinfo{booktitle}{Proceedings of the 2021 ACM conference on
  fairness, accountability, and transparency}, pp. \bibinfo{pages}{296--314}.
%Type = Article
\bibitem[{Ke et~al.(2017)Ke, Meng, Finley, Wang, Chen, Ma, Ye and
  Liu}]{ke2017lightgbm}
\bibinfo{author}{Ke, G.}, \bibinfo{author}{Meng, Q.}, \bibinfo{author}{Finley,
  T.}, \bibinfo{author}{Wang, T.}, \bibinfo{author}{Chen, W.},
  \bibinfo{author}{Ma, W.}, \bibinfo{author}{Ye, Q.}, \bibinfo{author}{Liu,
  T.Y.}, \bibinfo{year}{2017}.
\newblock \bibinfo{title}{Light{GBM}: A highly efficient gradient boosting
  decision tree}.
\newblock \bibinfo{journal}{Advances in Neural Information Processing Systems}
  \bibinfo{volume}{30}, \bibinfo{pages}{3149--3157}.
%Type = Article
\bibitem[{Kiefer(1953)}]{kiefer1953sequential}
\bibinfo{author}{Kiefer, J.}, \bibinfo{year}{1953}.
\newblock \bibinfo{title}{Sequential minimax search for a maximum}.
\newblock \bibinfo{journal}{Proceedings of the American Mathematical Society}
  \bibinfo{volume}{4}, \bibinfo{pages}{502--506}.
%Type = Article
\bibitem[{Rossi et~al.(1996)Rossi, McCulloch and Allenby}]{rossi1996value}
\bibinfo{author}{Rossi, P.E.}, \bibinfo{author}{McCulloch, R.E.},
  \bibinfo{author}{Allenby, G.M.}, \bibinfo{year}{1996}.
\newblock \bibinfo{title}{The value of purchase history data in target
  marketing}.
\newblock \bibinfo{journal}{Marketing Science} \bibinfo{volume}{15},
  \bibinfo{pages}{321--340}.
%Type = Inproceedings
\bibitem[{Subramanian et~al.(2022)Subramanian, Sun, Drissi and
  Ettl}]{subramanian2022constrained}
\bibinfo{author}{Subramanian, S.}, \bibinfo{author}{Sun, W.},
  \bibinfo{author}{Drissi, Y.}, \bibinfo{author}{Ettl, M.},
  \bibinfo{year}{2022}.
\newblock \bibinfo{title}{Constrained prescriptive trees via column
  generation}, in: \bibinfo{booktitle}{Proceedings of the AAAI Conference on
  Artificial Intelligence}, pp. \bibinfo{pages}{4602--4610}.
%Type = Article
\bibitem[{Thiele(2006)}]{thiele2006single}
\bibinfo{author}{Thiele, A.}, \bibinfo{year}{2006}.
\newblock \bibinfo{title}{Single-product pricing via robust optimization}.
\newblock \bibinfo{journal}{Submitted to Operations Research} .
%Type = Article
\bibitem[{Uehara et~al.(2024)Uehara, Nishimura, Li, Yang, Jobson, Ohashi,
  Matsumoto, Sukegawa and Takano}]{uehara2024robust}
\bibinfo{author}{Uehara, Y.}, \bibinfo{author}{Nishimura, N.},
  \bibinfo{author}{Li, Y.}, \bibinfo{author}{Yang, J.},
  \bibinfo{author}{Jobson, D.}, \bibinfo{author}{Ohashi, K.},
  \bibinfo{author}{Matsumoto, T.}, \bibinfo{author}{Sukegawa, N.},
  \bibinfo{author}{Takano, Y.}, \bibinfo{year}{2024}.
\newblock \bibinfo{title}{Robust portfolio optimization model for electronic
  coupon allocation}.
\newblock \bibinfo{journal}{arXiv preprint arXiv:2405.12865} .
%Type = Article
\bibitem[{Waldfogel(2015)}]{waldfogel2015first}
\bibinfo{author}{Waldfogel, J.}, \bibinfo{year}{2015}.
\newblock \bibinfo{title}{First degree price discrimination goes to school}.
\newblock \bibinfo{journal}{The Journal of Industrial Economics}
  \bibinfo{volume}{63}, \bibinfo{pages}{569--597}.
%Type = Article
\bibitem[{Wang et~al.(2021a)Wang, Chen and Wang}]{wang2021distributionally}
\bibinfo{author}{Wang, T.}, \bibinfo{author}{Chen, N.}, \bibinfo{author}{Wang,
  C.}, \bibinfo{year}{2021}a.
\newblock \bibinfo{title}{Distributionally robust prescriptive analytics with
  wasserstein distance}.
\newblock \bibinfo{journal}{arXiv preprint arXiv:2106.05724} .
%Type = Article
\bibitem[{Wang et~al.(2021b)Wang, Huang, Han and Lim}]{wang2021price}
\bibinfo{author}{Wang, X.}, \bibinfo{author}{Huang, H.C.},
  \bibinfo{author}{Han, L.}, \bibinfo{author}{Lim, A.}, \bibinfo{year}{2021}b.
\newblock \bibinfo{title}{Price optimization with practical constraints}.
\newblock \bibinfo{journal}{arXiv preprint arXiv:2104.09597} .
%Type = Article
\bibitem[{Yanagi et~al.(2024)Yanagi, Ikeda and Takano}]{yanagi2024robust}
\bibinfo{author}{Yanagi, T.}, \bibinfo{author}{Ikeda, S.},
  \bibinfo{author}{Takano, Y.}, \bibinfo{year}{2024}.
\newblock \bibinfo{title}{Robust portfolio optimization for recommender systems
  considering uncertainty of estimated statistics}.
\newblock \bibinfo{journal}{arXiv preprint arXiv:2406.10250} .

\end{thebibliography}




% \end{thebibliography}
\end{document}